\documentclass[a4paper,11pt]{article}

\usepackage[pdftex]{graphicx}
\usepackage{amsmath,amssymb}
\usepackage{amsthm}
\usepackage{subfig}
\usepackage[left=2.5cm, right=2.5cm]{geometry}
\usepackage[pagebackref]{hyperref}

\def\R{\mathbb{R}}

\begin{document}

\renewcommand{\thefootnote}{\fnsymbol{footnote}}

\begin{center}
{\Large \bf Travelling wave solutions of multisymplectic discretizations \\of semi-linear wave equations 
}

\vspace{3ex}
{\sc \small Fleur McDonald\footnotemark[1], Robert I McLachlan\footnotemark[2], 
Brian E Moore\footnotemark[3], and G R W Quispel\footnotemark[4]\footnotemark[5]}

\footnotetext[1]{Institute of Fundamental Sciences, Massey University, New Zealand
(\texttt{x\_999999999\_x@hotmail.com‎}).}

\footnotetext[2]{Institute of Fundamental Sciences, Massey University, New Zealand
($\mathtt{R.McLachlan@massey.ac.nz}$).
Work of this author was supported by the Marsden Fund of the Royal Society of New Zealand.
}

\footnotetext[3]{Department of Mathematics, University of Central Florida, USA
($\mathtt{Brian.Moore@ucf.edu}$).}

\footnotetext[4]{Department of Mathematics and Statistics, La Trobe University, Australia. 
($\mathtt{R.Quispel@latrobe.edu.au}$).
Work of this author was supported by the Australian Research Council.}

\
\renewcommand{\thefootnote}{\arabic{footnote}}

\end{center}

\vspace{1ex}
\noindent
{\bf Abstract.} How well do multisymplectic discretisations preserve travelling wave solutions?
To answer this question, the 5-point central difference scheme is applied to the semi-linear wave 
equation.  A travelling wave ansatz leads to an ordinary difference equation, 
whose solutions correspond to the numerical scheme and can be compared to travelling wave solutions 
of the corresponding PDE.  
For a discontinuous nonlinearity the difference equation is solved exactly.  
For continuous nonlinearities the difference equation is solved using a Fourier series, 
and resonances that depend on the grid-size are revealed for a smooth nonlinearity.
In general, the infinite dimensional functional equation, which must be solved to get the travelling wave solutions,
is intractable, but backward error analysis proves to be a powerful tool, as it provides a way to study the solutions 
of the equation through a simple ODE that describes the behavior to arbitrarily high order.
A general framework for using backward 
error analysis to analyze preservation of travelling waves for other equations and discretisations is presented.
Then, the advantages that multisymplectic methods have over other methods are briefly highlighted.  

\vspace{4ex}
\noindent
{\bf Key words.} five-point centered difference; backward error analysis; travelling wave solution; 
semi-linear wave equation; resonance 

\vspace{2ex}
\noindent
{\bf AMS subject classifications.} 65M22; 65P10; 35C07; 35L05 

\vspace{3ex}

\section{Introduction}
\label{MSIPDE}

Since their introduction by Marsden et al. \cite{30MGV} and Reich \cite{31MSR},
multisymplectic integrators have seen a huge growth in interest, because
they offer the prospect of getting the same excellent long-time results for Hamiltonian PDEs
that symplectic integrators do for ODEs \cite{32NMH,3SHD,5GIO}. The Hamiltonian ODE
$K z_t  = \nabla S(z)$, where $z\in\R^n$, $K$ is an invertible $n\times n$ antisymmetric
matrix and $S\colon\R^n\to \R$ is the Hamiltonian, obeys conservation of symplecticity:
$\omega_t = 0$, where $\omega = dz \wedge K dz$. The maps $z_n\to z_{n+1}$ obtained from
symplectic integrators obey $\Delta_t \omega = 0$. In contrast, the multi-Hamiltonian PDE
\begin{equation}
\label{eq:PDE}
K z_t + L z_x = \nabla S(z)
\end{equation}
where $K$ and $L$ are antisymmetric $n\times n$ matrices, obeys the multisymplectic conservation law
$\omega_t + \kappa_x = 0$,
where $\omega = dz\wedge K  dz$ and $\kappa = dz\wedge L dz$. Multisymplectic integrators
have a fully discrete conservation law which is a discretization of the continuous one.

There are a number of apparent advantages to multisymplectic integrators.
They are in some sense the simplest natural discrete analogue of (\ref{eq:PDE}).
Many standard methods that have been in use for a long time are multisymplectic,
such as the 5-point stencil for the wave equation. They obey a discrete variational principle.
Because the conservation law is local, they allow a flexible treatment of boundary conditions.
They are nondissipative. They appear to work well in many numerical experiments on integrable
and nonintegrable PDEs.

However, some differences between the PDE and ODE case are immediately apparent.
(i) In the PDE case, the matrices $K$ and $L$ can be singular,
and so the 2-forms $\omega$ and $\kappa$ can be degenerate.
(ii) In the ODE case,
the discretization of the conservation law is unique, and is identical to the one obeyed by
the flow of the ODE, whereas in the PDE case, the discretization of the conservation law
is not unique and is necessarily an approximation of that obeyed by the PDE.
(iii) In the ODE,
preservation of the symplectic form is known to be fundamental to the dynamics of Hamiltonian systems,
and is used either implicitly or explicitly in all studies of Hamiltonian dynamics. In the PDE,
use of a multisymplectic conservation law is hardly known. It is not clear if this is because the
study of the dynamics of PDEs itself is relatively less developed, or because the conservation
law itself is relatively weak. The results and observations available for multisymplectic
integrators so far are not as strong as those available for symplectic integrators, which
include backward error analysis, long-time energy behaviour, and preservation of
invariant sets such as periodic and quasiperiodic orbits.

Nevertheless, multisymplectic schemes have been shown to preserve more than just the
multisymplectic conservation law.  For example, the Gauss--Legendre Runge--Kutta methods
satisfy fully discrete energy and momentum conservation laws for linear equations \cite{15MSI,31MSR}.
The Preissmann box scheme and the Euler box scheme satisfy
semi-discrete energy and momentum conservation laws for nonlinear equations \cite{63BEA,41MSI}.
Other properties preserved by certain multisymplectic methods include phase space structure \cite{IS04},
conservation of wave action \cite{JF06}, dispersion relations and group velocity \cite{45MBS,209LPN},
and other conservation laws that result from Noether's theorem \cite{CHH07,41MSI}.

In this paper, we study the preservation of travelling waves under multisymplectic discretization.
Many PDEs, such as the Boussinesq, Schr\"{o}dinger, Korteweg--de Vries, and
nonlinear wave equations (among others), have travelling wave solutions,
and the multi-Hamiltonian forms for these conservative PDEs are known (cf.~\cite{TB1,TB2,BD99,BD01}).
(Note that some equations can be written in the multi-Hamiltonian form (\ref{eq:PDE}) in multiple ways \cite{77MSR}.)

Our focus is on the semi-linear wave equation
\begin{equation}
\label{eq:NLW}
u_{tt} = u_{xx} -V'(u).
\end{equation}
Travelling waves are solutions to the PDE that travel at a constant speed $c$ without changing shape.
Thus, they take the form $u(x,t) = \varphi(x-ct)$. The PDE \eqref{eq:NLW} has a 2-parameter group of
symmetries given by translations in $x$ and $t$; thus travelling waves are relative equilibria with
respect to a subgroup of this symmetry group. Introducing the travelling wave coordinate
$\xi := x-ct$, the PDE (\ref{eq:NLW}) is reduced to the ODE
\begin{equation}
\label{eq:reduced}
(c^2-1)\varphi_{\xi\xi}(\xi) = -V'(\varphi(\xi)).
\end{equation}
This is a planar Hamiltonian system. Introducing the conjugate momentum $\psi := (c^2-1)\varphi_\xi$, with
$(c^2-1)$ playing the role of mass, its Hamiltonian form is
\begin{equation}
\label{eq:reducedH}
\varphi_\xi = \frac{\partial H}{\partial\psi}, \qquad
\psi_\xi =-\frac{\partial H}{\partial\varphi} \qquad \textup{with} \qquad
H = \frac{1}{2}(c^2-1)^{-1}\psi^2 + V(\varphi).
\end{equation}
Its phase portrait is easily determined from the level sets of the reduced Hamiltonian $H$.
Bounded solutions are either
solitary waves (that obey $\varphi'(\pm\infty)=0$) or periodic travelling waves.
Solitary waves (sometimes called pulses, wave fronts, kinks, or anti-kinks) correspond to homoclinic
($\varphi(\infty)=\varphi(-\infty)$) or heteroclinic ($\varphi(\infty)\ne\varphi(-\infty)$) orbits of
(\ref{eq:reduced}), and require an infinite spatial domain. Periodic travelling waves correspond to
periodic solutions of (\ref{eq:reduced}) and may exist on a periodic or infinite spatial domain.

Multisymplectic integrators for equation (\ref{eq:NLW}) have been constructed using discrete variational principles
\cite{30MGV}, and by applying symplectic Runge--Kutta \cite{MRS,31MSR,29MSI}
and partitioned Runge--Kutta \cite{208MSP,206LSP,29MSI,33OMP} methods in space and time.
In this paper we consider one particular multisymplectic integrator,
the 5-point central difference method for (\ref{eq:NLW}) given by
\begin{equation}
\label{eq:5pt}
\frac{1}{(\Delta t)^2}\left(u_{i}^{n+1}-2u_{i}^{n}+u_{i}^{n-1}\right) -
\frac{1}{(\Delta x)^2}\left(u_{i+1}^{n}-2u_{i}^{n}+u_{i-1}^{n}\right) = - V'(u_{i}^{n}).
\end{equation}
This can be viewed either as the 2-stage Lobatto IIIA--IIIB partitioned Runge--Kutta method
applied in space and time with a particular choice of partitioning \cite{33OMP},
as the Euler-Lagrange equations of a discrete Lagrangian \cite{63BEA},
or simply as the leapfrog method in space and time.

We seek  travelling wave solutions of the discrete equation (\ref{eq:5pt}) of the form
\begin{equation}
u_i^n = \varphi(x_i - c t_n) = \varphi(i \Delta x - c n \Delta t)= \varphi(i\sigma-n\kappa)= \varphi(\xi)
\end{equation}
where $\sigma = \Delta x$, $\kappa = c \Delta t$, and $\xi = i \sigma - n \kappa$.
The two parameters $\sigma$ and $\kappa$ will play an important role in the description of discrete travelling waves.
The substitution $u_i^n=\varphi(\xi)$ into the scheme (\ref{eq:5pt}) yields the
 {\em discrete travelling wave equation}
\begin{align} \label{DTWE}
& \frac{c^{2}}{\kappa^{2}}\left(\varphi(\xi+\kappa)-2\varphi(\xi)+\varphi(\xi-\kappa)\right) \\ \notag
&\qquad\qquad -\frac{1}{\sigma^{2}}\left(\varphi(\xi+\sigma)-2\varphi(\xi)+\varphi(\xi-\sigma)\right) =-V'(\varphi(\xi)).
\end{align}
Ideally, one would like to understand all solutions of (\ref{DTWE}) under broad conditions on
$c$, $\kappa$, $\sigma$, $V$, and $\varphi$---for example, $V$ smooth and $\varphi$
bounded for all $\xi\in \R$. This entire family of discrete travelling waves could then be
compared to the continuous family. It is immediately apparent that the parameter $\sigma/\kappa$
is of crucial importance. When $\sigma/\kappa=m/n$ is rational, (\ref{DTWE}) is a finite-dimensional map.
It then makes the most sense to look for solutions with $\xi\in \frac{\sigma}{m}\mathbb{Z}$.
When $\sigma/\kappa$ is irrational, we need to take $\xi\in\mathbb{R}$ and seek solutions in some
space of functions to be determined from (\ref{DTWE}) itself. In this case, (\ref{DTWE}) is an infinite-dimensional
functional equation. Thus, at first sight, the drastic reduction in complexity that is obtained in the
continuous case (from a PDE to a planar ODE) is not obtained in the discrete case.

We approach equation (\ref{DTWE}) from two distinct points of view; hence, the article is divided into two parts.
In Part I, {\sl case studies} are used to construct solutions of (\ref{DTWE}) and 
carefully compare them to the solutions of (\ref{eq:reduced})
for three special cases of the nonlinearity: \vspace{-0.25ex}
\begin{enumerate}
\item a discontinuous $V^{\prime}$ given by the McKean caricature of the cubic, in Section \ref{sec:mckean}, \vspace{-1ex}
\item a continuous but non-smooth $V^{\prime}$ given by the sawtooth nonlinearity, in Section \ref{sec:sawtooth}, \vspace{-1ex}
\item a smooth $V^{\prime}$ given by the sine nonlinearity, in Section \ref{sec:smooth}.  \vspace{-0.25ex}
\end{enumerate}
Our main contributions in Part I are:
\begin{itemize}
\item Comparison of discrete and continuous travelling waves providing details of the extent to which
travelling waves are preserved by the discretisation. \vspace{-0.5ex}
\item Exposition of tools and techniques for constructing and analyzing solutions of (\ref{DTWE}). \vspace{-0.5ex}
\item Exploration of resonances that depend on $\kappa$ and $\sigma$ for a smooth nonlinearity. \vspace{-0.5ex}
\end{itemize}
In Part II, {\sl backward error analysis} is used to give a highly accurate description of the discrete travelling waves, 
as well as the way in which they approximate the continuous travelling waves.  
This more general approach allows for considerations of a wider range of problems, which is illustrated in two sections.
\begin{enumerate}
 \item Application of backward error analysis for (\ref{DTWE}) with arbitrary $V^{\prime}$, in Section \ref{sec:multistep}.
 \item Demonstration of how backward error analysis may be used to study travelling waves 
for other multisymplectic discretizations and multi-Hamiltonian equations, in Section \ref{sec:conclude}.
\end{enumerate}
Our main contributions in Part II are:
\begin{itemize}
\item Introduction of an original and extremely powerful application of backward error analysis to the study of discrete travelling waves, 
which reduces (in the sense of backward error) an intractable infinite-dimensional functional equation to a simple ODE.\vspace{-0.5ex}
\item Confirmation that the travelling wave solutions of the modified equations (obtained through backward error analysis)
mimic the travelling wave solutions of the PDE and represent a solution of the  multisymplectic discretisation.\vspace{-0.5ex}
\item Development of a framework for analyzing the travelling wave solutions of the modified 
equations obtained through backward error analysis for other multi-symplectic discretisations 
of more general multi-Hamiltonian PDEs. \vspace{-0.5ex}
\end{itemize}
Finally, a juxtaposition of multisymplectic and non-symplectic (especially symmetric) methods 
for the steady-state solutions (in the Appendix)
provides a useful and interesting contrast of the numerical solution behavior, revealing a special advantage 
of multisymplectic methods. 

\part{Case Studies}
Travelling wave solutions of the multi-symplectic discretization (\ref{eq:5pt}) are obtained 
by solving the ordinary difference equation (\ref{DTWE}) for some specific cases of $V^{\prime}$.
Similar studies have been performed in different contexts.  
Most relevant are the results of Cahn et al.~\cite{CMPVV}, which construct exact travelling wave 
solutions for a spatially discrete Allen-Cahn equation; since their equation has two space dimensions,
their steady-state equation is actually a special case of (\ref{DTWE}).  In addition, the equation solved 
in \cite{CMPVV} has a McKean nonlinearity, which is used in our first case study, and 
our method of solution in that case follows the approach of \cite{CMPVV}.

\section{Discontinuous $V^{\prime}$: McKean nonlinearity} \label{sec:mckean}

\begin{figure}[h]
\centering
\includegraphics[width=0.4\textwidth]{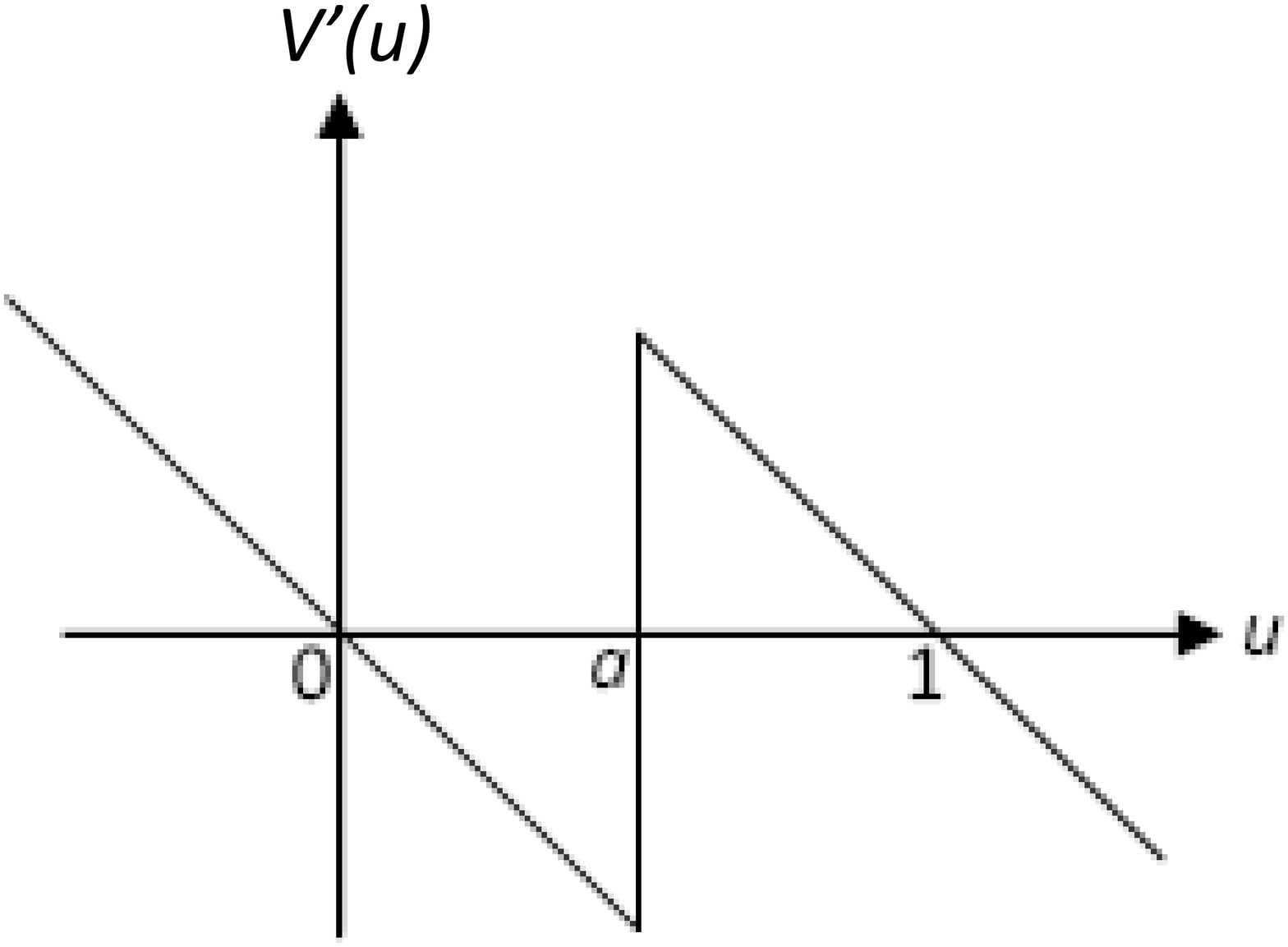}
\caption{McKean nonlinearity: $-V'(u) = -u + h(u-a)$.  \label{McKeanGraph}}
\end{figure}

In a 1970 study of nerve conduction, McKean \cite{McK} proposed a caricature of the cubic function $x(x-a)(x-1)$.
The McKean caricature (Figure ~\ref{McKeanGraph}) is given by
\begin{equation}
\label{McKCub}
V^{\prime}(u)=u-h(u-a)
\end{equation}
where $0<a<1$ and $h$ is the Heaviside step function.
The function is discontinuous, but it has long been used in the study of travelling waves,
and it provides a straightforward path to exact solutions of (\ref{eq:reduced}) and (\ref{DTWE}).

\subsection{Continuous travelling waves}
\label{McKeanCont}
\begin{figure}
\centering
\includegraphics[angle=270,width=0.6\textwidth]{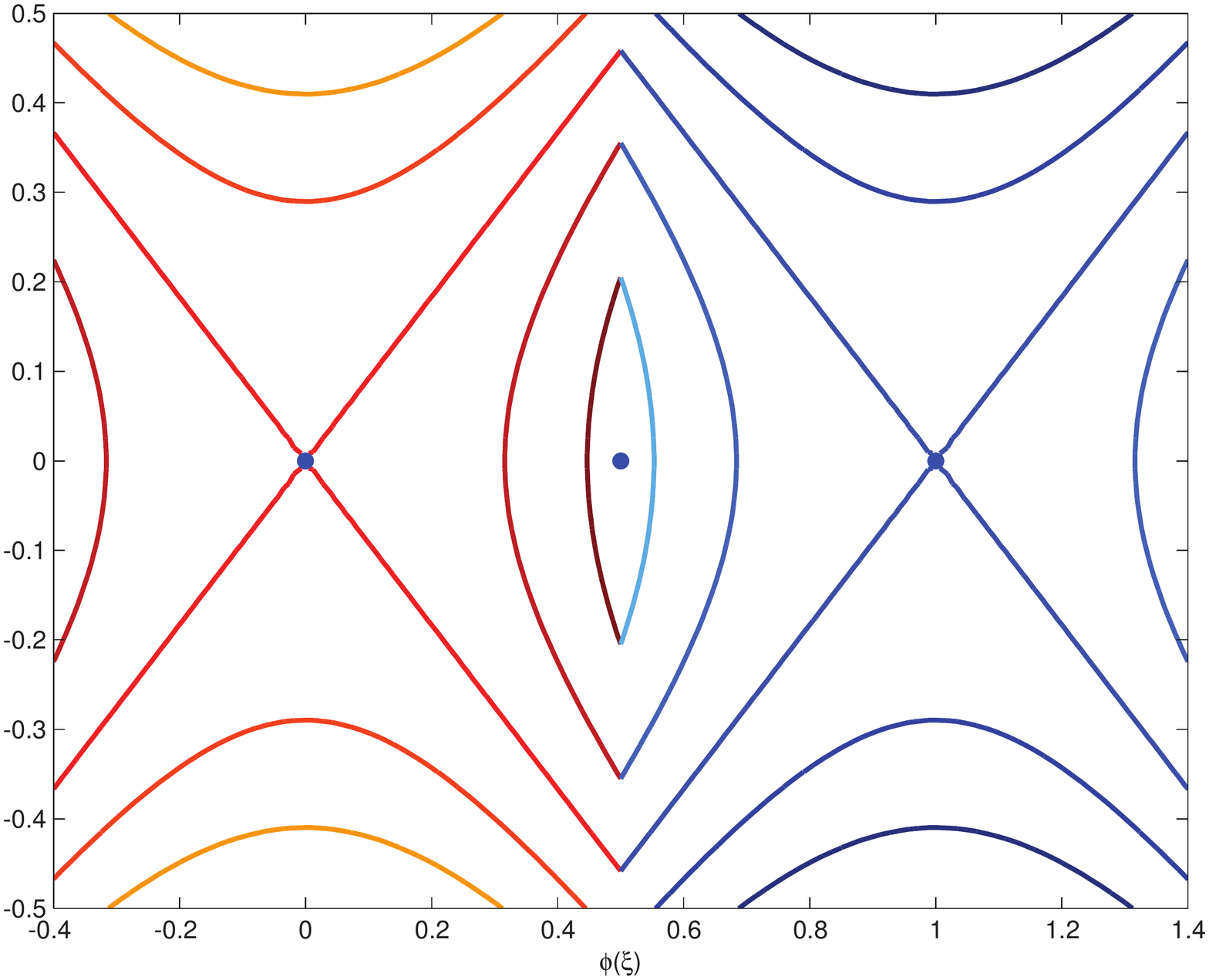}
\caption{Phase portrait of (\ref{eq:reduced}) with the McKean nonlinearity.  \label{McKeanPhase}}
\end{figure}

\noindent Upon substituting \eqref{McKCub} into \eqref{eq:reduced}, the travelling wave equation becomes
\begin{equation}
(c^{2} - 1)\varphi''(\xi) = -\varphi(\xi) + h(\varphi(\xi) - a),
\label{ContDD1}
\end{equation}
where we require $\left|c\right|<1$.
The phase portrait of this ODE is shown in Figure~\ref{McKeanPhase}.
Notice that \eqref{ContDD1} possesses both periodic and heteroclinic travelling waves.
For heteroclinic travelling waves we impose the boundary conditions
\begin{equation}
\lim_{\xi\to -\infty}\varphi(\xi) = 0, \qquad \lim_{\xi\to \infty}\varphi(\xi) = 1.
\label{BCs}
\end{equation}
We expect a monotonic front solution, and, since $a\in(0,1)$,
there must be $\xi^{*} \in \mathbb{R}$ for which
\begin{equation}
\label{hChange}
\varphi(\xi^{*}) = a \quad \text{with} \qquad 
\varphi(\xi)<a \quad \text{for} \quad \xi < \xi^{*} \qquad \textup{and} \qquad
\varphi(\xi)>a \quad \text{for} \quad \xi > \xi^{*}.
\end{equation}
Following \cite{CMPVV}, the continuous nonlinear problem (\ref{ContDD1}) may be stated,
\begin{equation}
(c^{2} - 1)\varphi''(\xi) = -\varphi(\xi) + h(\xi-\xi^{*}),
\label{ContDD2}
\end{equation}
which is linear.
This is the essential reason for considering the McKean nonlinearity.

We seek solutions of equation (\ref{ContDD2})
for which there exists an $\epsilon>0$ such that
$\left|\varphi(\xi)\right|\leq Ke^{\epsilon\xi}$ for $\xi<0$ (see Lemma 4.1 of \cite{CMPVV}). We apply the change of variables
\begin{equation}
\varphi_{\epsilon}(\xi) = e^{-\epsilon\xi}\varphi(\xi)
\label{COV}
\end{equation}
where $\epsilon>0$ is sufficiently small, and obtain
\begin{equation*}
(c^{2}-1)(2\epsilon\varphi'_{\epsilon}(\xi) + \varphi''_{\epsilon}(\xi)) + \varphi_{\epsilon}(\xi)(1+\epsilon^{2}(c^{2}-1)) =  e^{-\epsilon\xi}h(\xi-\xi^{*}).
\end{equation*}

Applying the Fourier transform,
\begin{equation}
\hat{\varphi}_{\epsilon}(s) = \int^{\infty}_{-\infty}e^{-is\xi}\varphi_{\epsilon}(\xi)d\xi
\label{FT}
\end{equation}
and using the properties of Fourier transforms, we get
\begin{equation*}
\hat{\varphi}_{\epsilon}(s) = \frac{e^{-(is+\epsilon)\xi^{*}}}{(is+\epsilon)R_{\text{Cont}}(s-i\epsilon)}
\end{equation*}
where $R_{\text{Cont}}(s) = 1 - c^{2}s^{2} + s^{2}$.
Next, we take the inverse Fourier transform,
\begin{equation}
\varphi_{\epsilon}(\xi) = \frac{1}{2\pi}\int^{\infty}_{-\infty}e^{is\xi}\hat{\varphi}_{\epsilon}(s)ds,
\label{IFT}
\end{equation}
giving the solution in the original variables as
\begin{align*}
\varphi(\xi)
&= \frac{1}{2\pi}\int^{\infty}_{-\infty}\frac{e^{(is+\epsilon)\xi}e^{-(is+\epsilon)\xi^{*}}}{(is+\epsilon)R_{\text{Cont}}(s-i\epsilon)}ds.
\end{align*}

\begin{figure}
\centering
\includegraphics[angle = 270,width=0.3\textwidth]{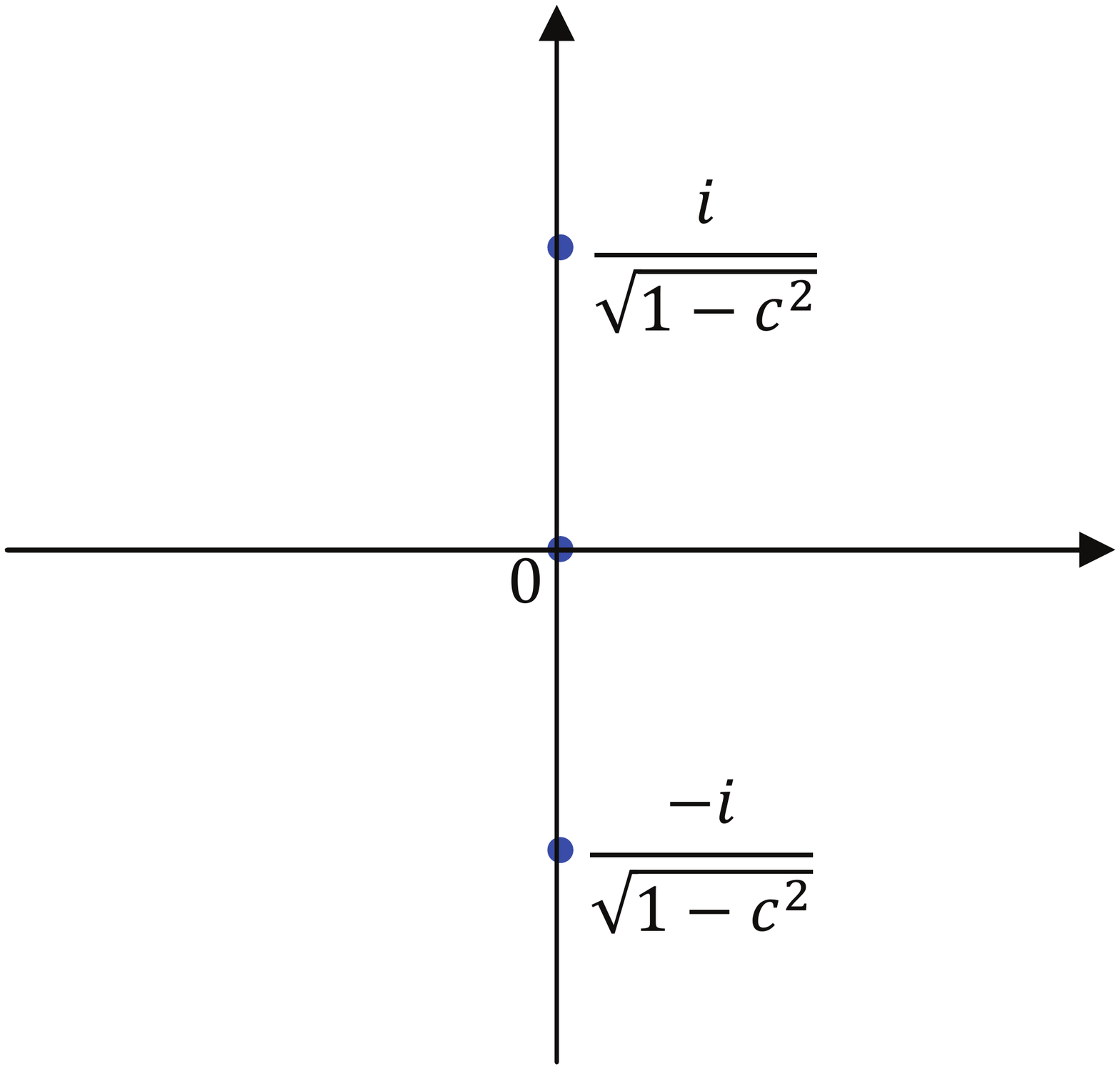}
\hspace{10mm}
\includegraphics[angle = 270,width=0.3\textwidth]{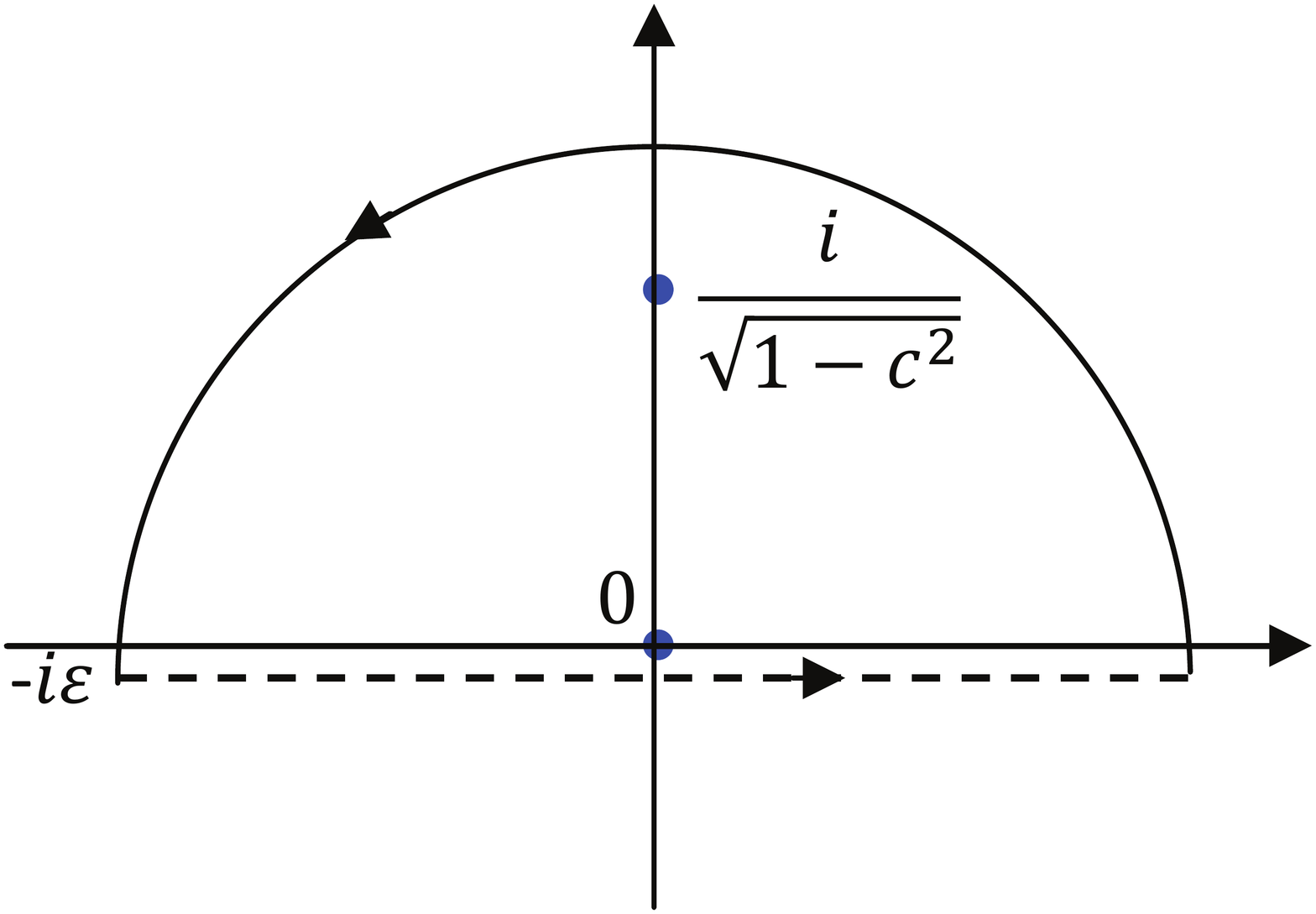}
\vspace{-10pt}
\caption{Singularities in the integrand of \eqref{Fintegral} and the contour used. 	\label{UpperContour}}
\end{figure}

Taking the limit as $\epsilon\to 0$ and changing the limits of integration we get
\begin{align} \label{Fintegral}
\varphi(\xi)
&= \frac{1}{2\pi i}\int^{-i\epsilon + \infty}_{-i\epsilon - \infty}\frac{e^{is(\xi-\xi^{*})}}{sR_{\text{Cont}}(s)}ds.
\end{align}
The domain of integration is changed slightly in order that residue theory can be applied to evaluate the improper integral.
We can factorise $R_{\text{Cont}}(s)$ to get
\begin{equation*}
R_{\text{Cont}}(s) = \left(1-c^{2}\right)\left(s-\frac{i}{\sqrt{1-c^{2}}}\right)\left(s+\frac{i}{\sqrt{1-c^{2}}}\right).
\end{equation*}
Now, we evaluate the integral, noting singularities
\begin{equation*}
s_{0} = 0, \qquad s_{1} = \frac{i}{\sqrt{1-c^{2}}}, \qquad s_{2} = -\frac{i}{\sqrt{1-c^{2}}},
\end{equation*}
which are plotted in Figure~\ref{UpperContour} for $0\leq c<1$.

\begin{figure}[h]
\centering
\includegraphics[angle=270,width=0.55\textwidth]{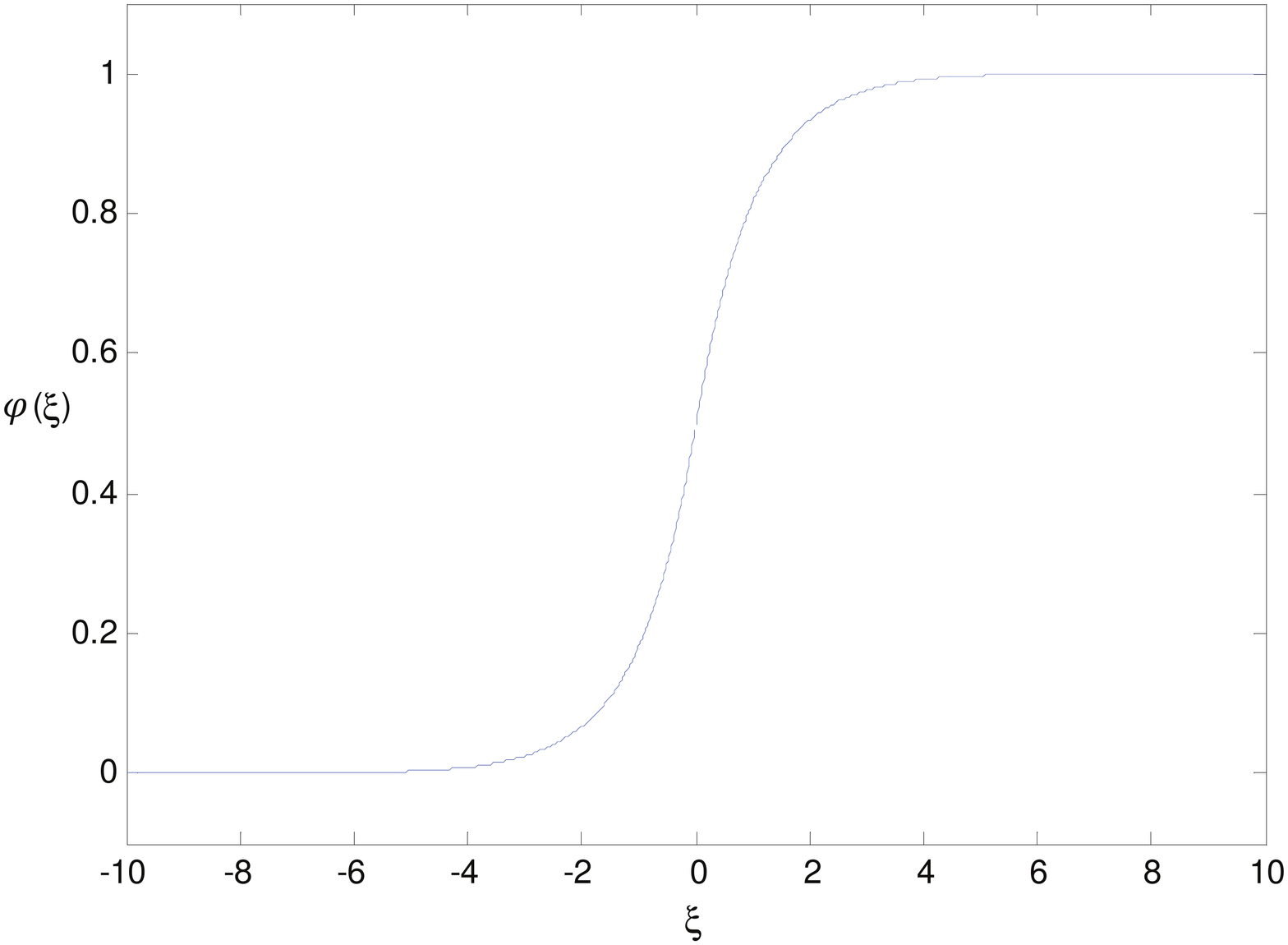}
\caption{Travelling wave solution \eqref{ContSolutionEq} of \eqref{eq:reduced} with $c=0.1$ and $\xi^{*}=0$.  \label{ContSolution}}
\end{figure}

For $\xi>\xi^*$ we close the top half plane with a semi-circular region taking the real axis to $-i\epsilon$.
(See Figure~\ref{UpperContour}.) The singularities enclosed by this region are $s_{0}$ and $s_{1}$.
For $\xi<\xi^*$ we close the bottom half plane enclosing $s_2$.
Hence, the explicit analytic solution for the continuous case is
\begin{equation}
\label{ContSolutionEq}
\varphi(\xi) = \left\{
\begin{array}{rl}
1 - \frac{1}{2}e^{\frac{-\left(\xi-\xi^{*}\right)}{\sqrt{1-c^{2}}}} & \text{if } \xi > \xi^{*}\\
\frac{1}{2}e^{\frac{\left(\xi-\xi^{*}\right)}{\sqrt{1-c^{2}}}} & \text{if } \xi < \xi^{*}
\end{array} \right.
\end{equation}
which is a heteroclinic travelling wave, shown in  Figure~\ref{ContSolution} for $c=0.1$ and $\xi^{*}=0$.

\subsection{Discrete travelling waves}

Now consider the discrete travelling wave equation \eqref{DTWE} with McKean nonlinearity.
Using the same procedure as  for the continuous case gives the solution
\begin{equation}
\label{SolMcKeanDis}
\varphi(\xi) = \frac{1}{2\pi i}\int^{\infty}_{-\infty}\frac{e^{is(\xi-\xi^{*})}}{sR_{\text{Disc}}(s)}ds
\end{equation}
where
\begin{equation}
R_{\text{Disc}}(s) = 1 - \frac{2c^{2}}{\kappa^{2}}\left(1 - \cos(\kappa s)\right) + \frac{2}{\sigma^{2}}\left(1 - \cos(\sigma s)\right).
\label{RDiscrete}
\end{equation}

Finding the explicit solution $\varphi(\xi)$ for this discrete case comes down to
analyzing $R_{\text{Disc}}(s)$. For rational values of $\frac{\sigma}{\kappa}$,
finding the zeros of $R_{\text{Disc}}(s)$ reduces to finding the roots of polynomials.
We illustrate this in two cases.

\begin{figure}[h]
\centering
\includegraphics[width=0.49\textwidth]{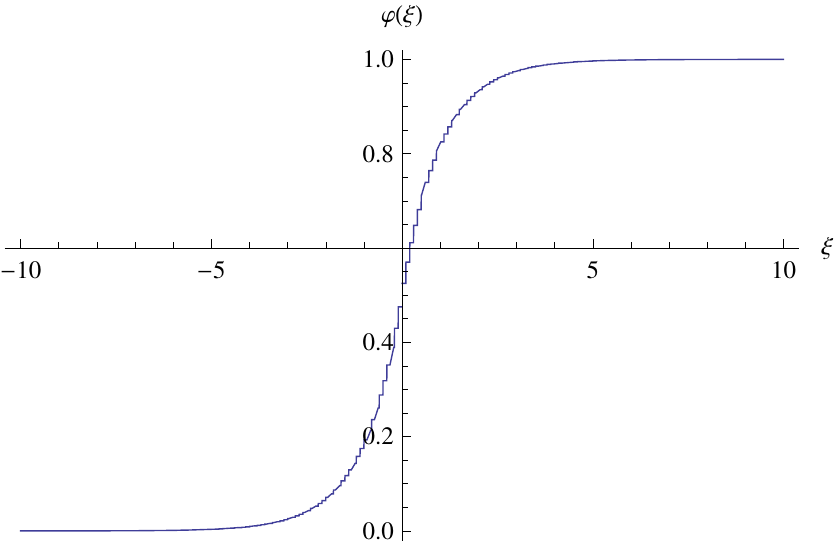}
\includegraphics[width=0.49\textwidth]{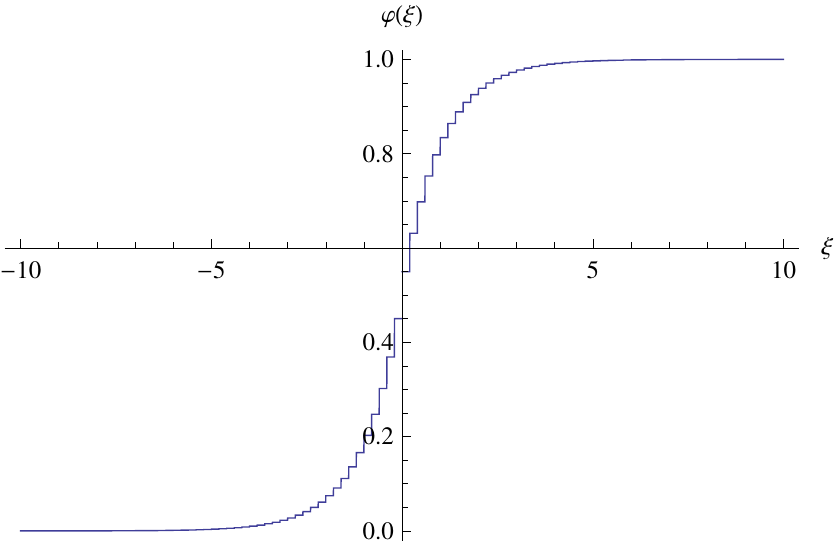}
\caption{The solution of \eqref{DTWE} for $c=0$: $\kappa=0.1$ left, $\kappa=0.2$ right.  \label{McKeanc0}}
\end{figure}

\begin{figure}[h]
\centering
\includegraphics[width=0.49\textwidth]{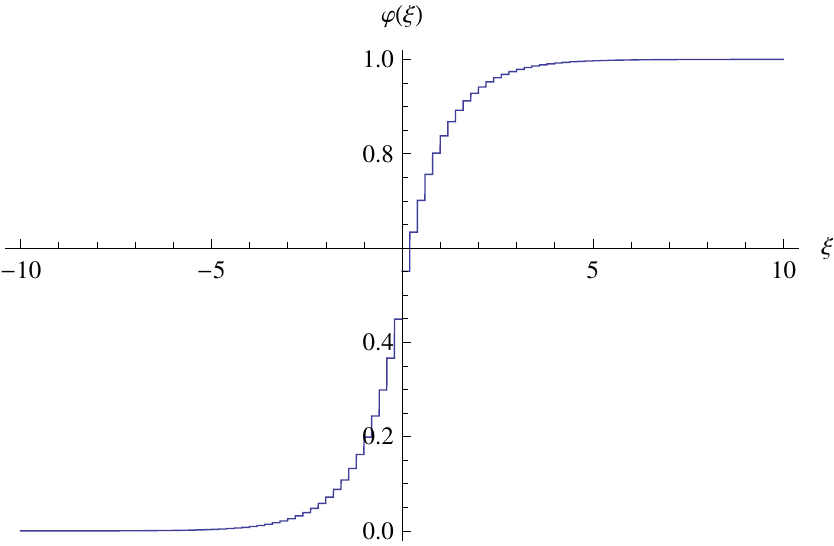}
\includegraphics[width=0.49\textwidth]{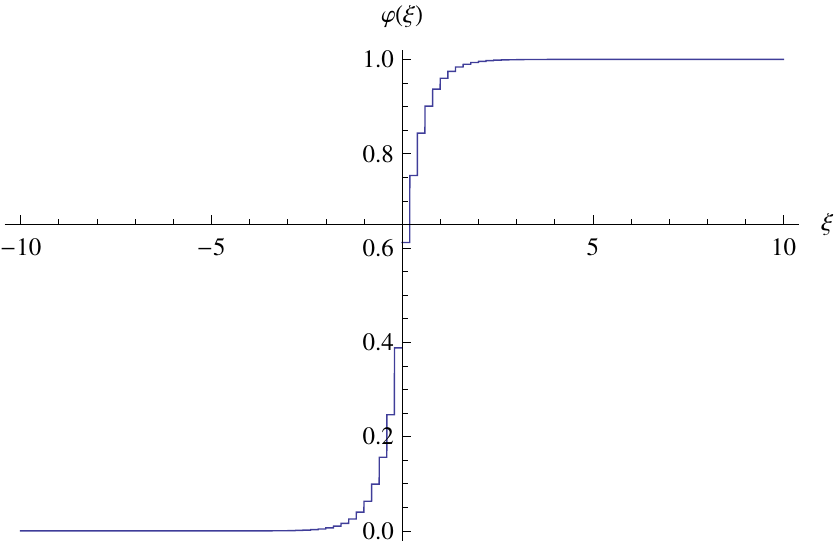}
\caption{The solution of \eqref{DTWE} for $\sigma=\kappa$: $c=0.2$ left, $c=0.9$ right.  \label{McKeansk}}
\end{figure}

\noindent {\bf Case $\sigma=\kappa$.}
Setting $\sigma=\kappa$ in \eqref{RDiscrete} we get
\begin{equation*}
R_{\text{Disc}}(s)=1+\frac{2}{\kappa^2}\left(1-\cos(\kappa s)\right)\left(1-c^2\right)
\end{equation*}
which has zeros at 
\begin{equation*}
s=\frac{\pm\arccos\left(1+\frac{\kappa^{2}}{2(1-c^2)}\right)\pm 2n\pi}{\kappa}
\end{equation*}
where $n=1,2,3,\cdots$.
The residues can be calculated and summed analytically, giving an explicit solution of the discrete travelling wave equation.
Examples for $c=0$ and for $c\neq 0$ are shown in Figures \ref{McKeanc0} and \ref{McKeansk}, respectively.
In all cases the solution is monotonic and piecewise constant, and the tail at each end of the solution does not contain any wiggles.

\begin{figure}[h]
\centering
\includegraphics[angle=270,width=0.55\textwidth]{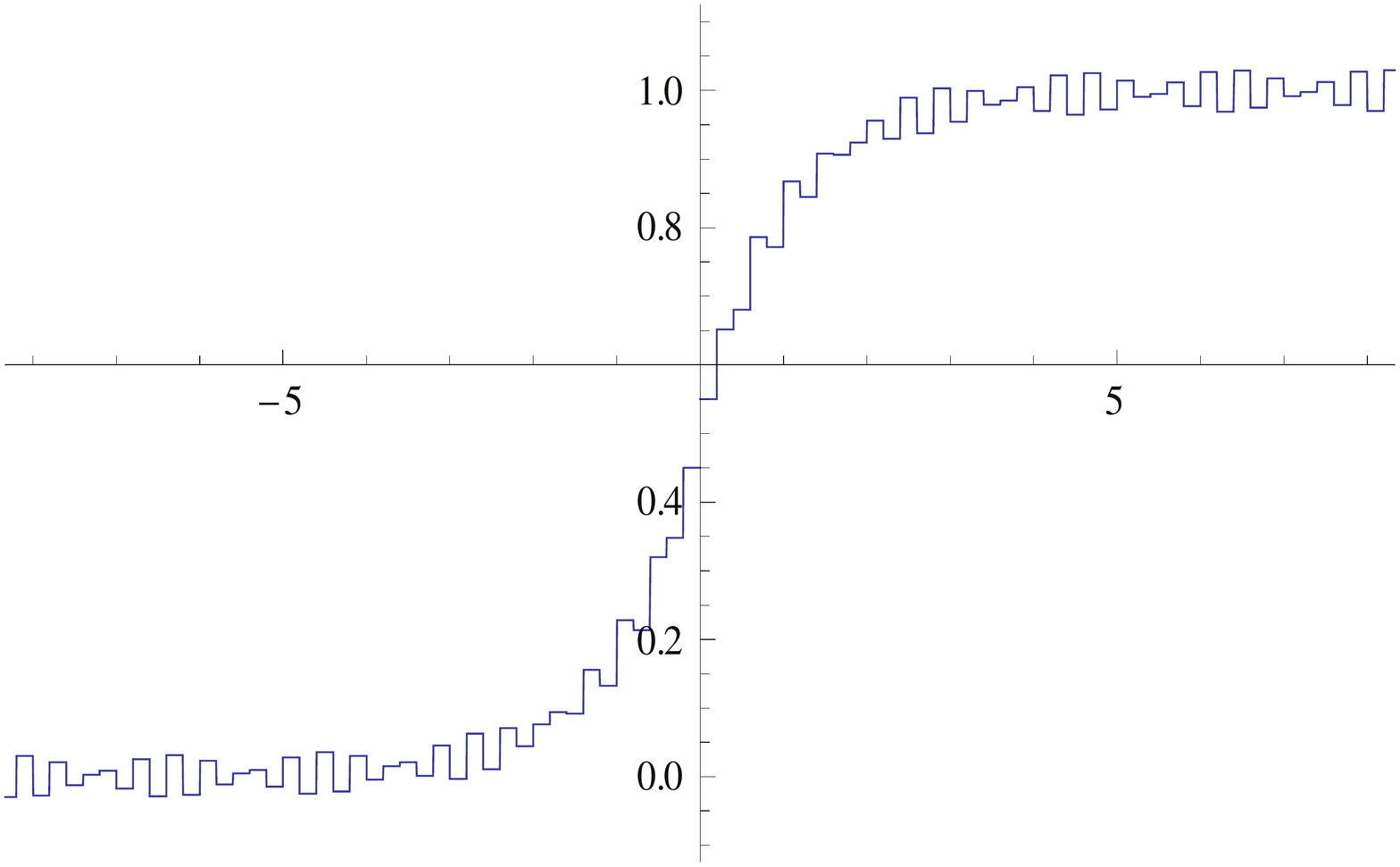}
\vspace{-20pt}
\caption{Solution of the discrete travelling wave equation with McKean nonlinearity for $\sigma=2\kappa$.  \label{DiscreteWave}}
\end{figure}

\noindent {\bf Case $\sigma=2\kappa$.}
We now get
\begin{equation*}
R_{\text{Disc}}(s) = 1 - \frac{2c^{2}}{\kappa^{2}}\left(1 - \cos(\kappa s)\right) + \frac{1}{2\kappa^{2}}\left(1 - \cos(2\kappa s)\right).
\end{equation*}
The zeros of $R_{\text{Disc}}(s)$ can be found explicitly.
There are  two complex roots and two real roots, plus their periodic images.
The sum of the residues over the periodic images can be calculated explicitly in terms of
hypergeometric functions. A typical solution is shown in Figure~\ref{DiscreteWave}.
Notice that we get a piecewise constant solution with wiggles at the tails.
These wiggles, originating from the real zeros of $R_{\text{Disc}}(s)$,
have amplitude $\mathcal{O}(\kappa^2)$.

In theory, such an analytic solution can be found for any rational value of $\sigma/\kappa$,
and in principle one could take a sequence of rational values approaching an irrational,
although this would be very complicated. It is difficult to make conclusions about the
existence of discrete travelling waves for irrational $\sigma/\kappa$ directly because this
requires analyzing the convergence of the integral \eqref{SolMcKeanDis},
which in turn requires detailed information about the zeros of the quasiperiodic function $R_{\text{Disc}}(s)$.

\subsection{Periodic travelling waves}
\label{PeriodicSolutions}

Discrete periodic travelling waves of period $2\tau$, say, can be sought as Fourier series.
First, set $\xi\in[0,2\tau)$ so that
\begin{equation*}
h(\xi) = \left\{
\begin{array}{rl}
0 & \text{if } \xi \in [0,\tau)\\
1 & \text{if } \xi \in (\tau,2\tau).
\end{array} \right.
\end{equation*}

\begin{figure}
\begin{center}
\includegraphics[width=0.8\textwidth]{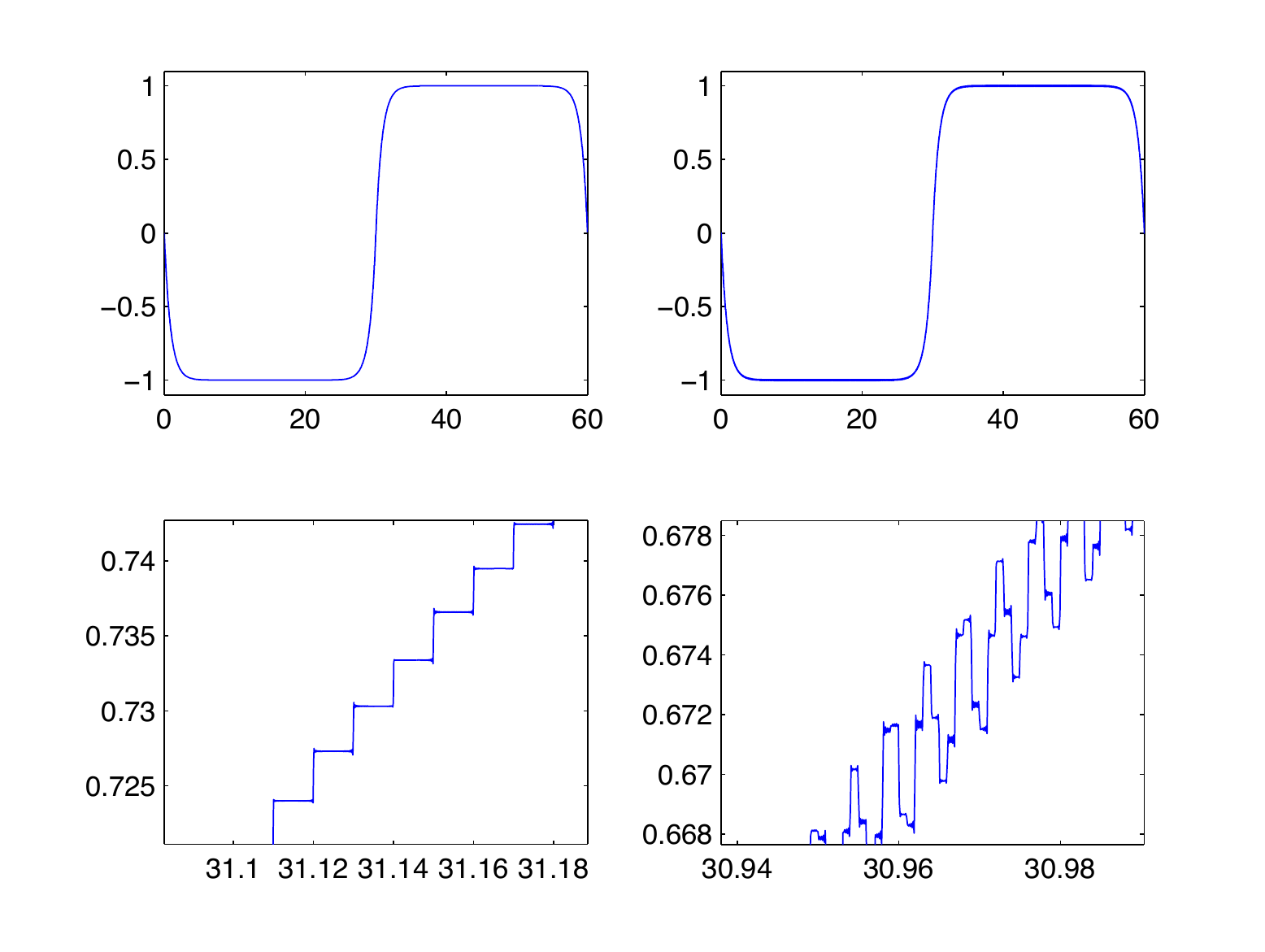}
\includegraphics[width=0.8\textwidth]{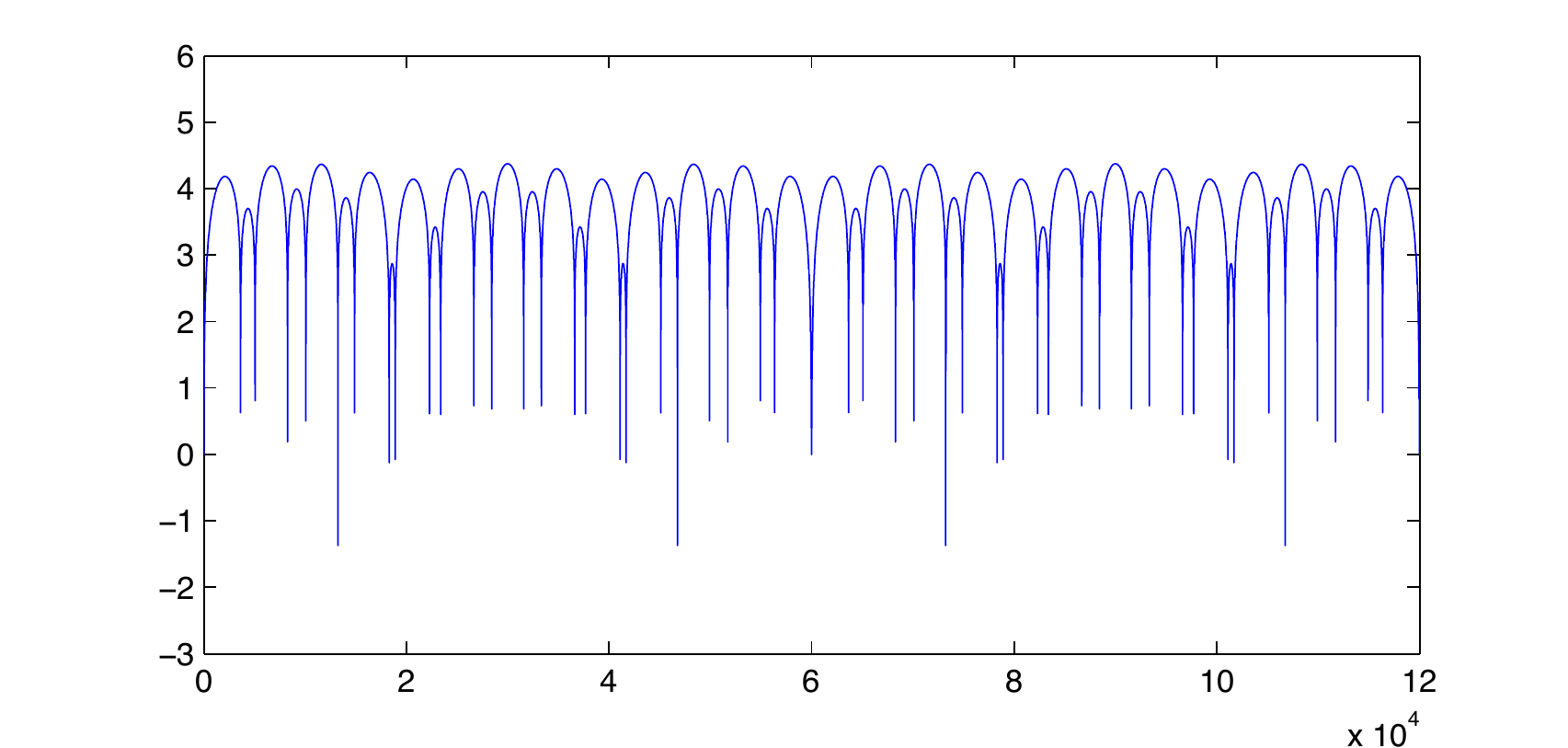}
\caption{\label{fig:rrat} Periodic travelling waves of the leapfrog method with the McKean nonlinearity.
When the space step $\sigma$, time step $\kappa$, and period $\tau$ are rationally related, travelling waves persist.
Here we take $\tau = 30$ and $c=1/2$ and approximate the waves with $2^{19}$ terms of a Fourier series.
The left hand column shows the solution for $\sigma/\kappa=2$, the right hand side for $\sigma/\kappa = 10/7$.
The bottom figure shows the base-10 logarithm of the Fourier amplitudes for $r=10/7$.
}
\end{center}
\end{figure}

\begin{figure}[h]
\begin{center}
\includegraphics[width=0.5\textwidth]{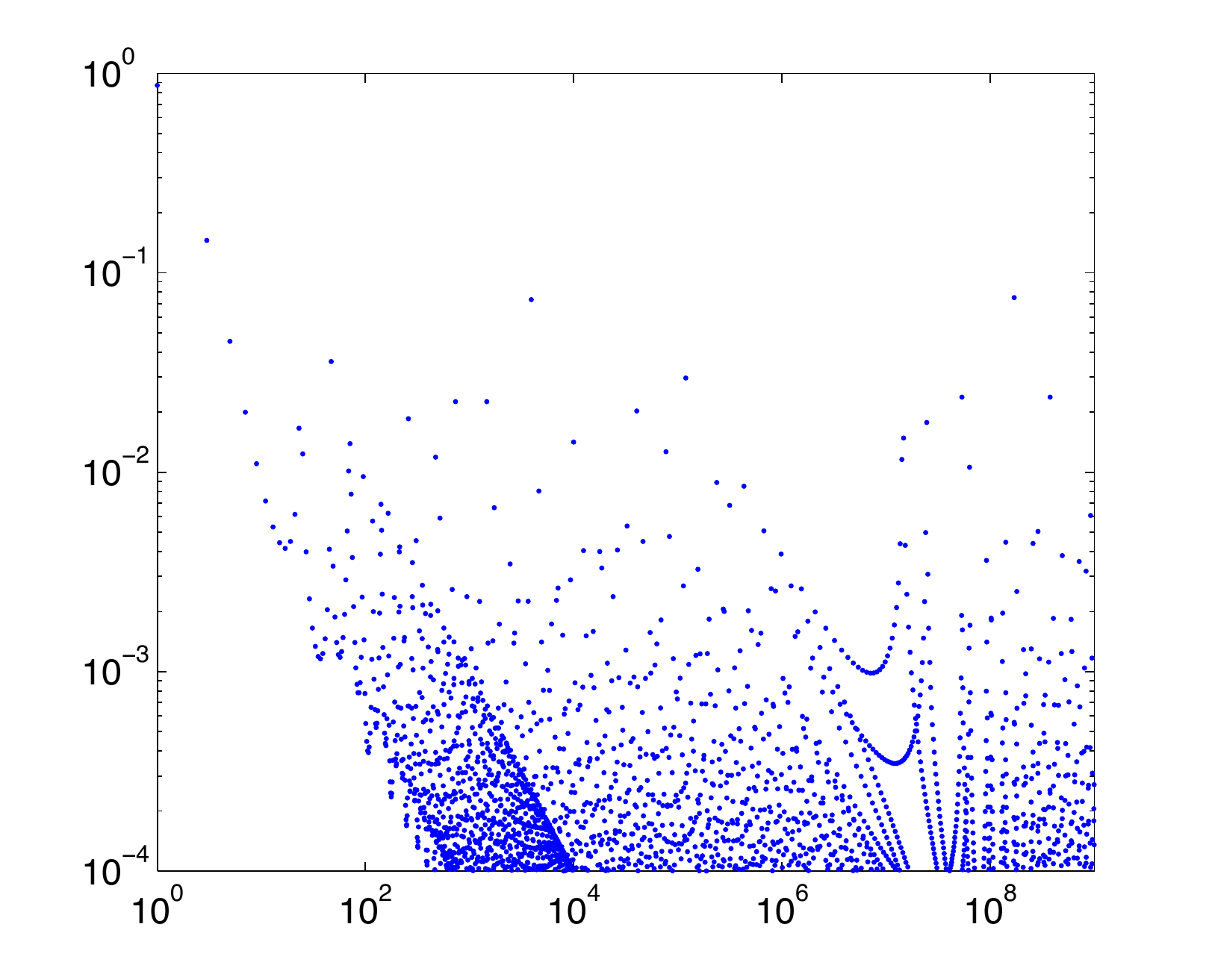}
\caption{\label{fig:PhiNonConverge}
Fourier amplitudes of a putative discrete travelling wave using a McKean nonlinearity with
$c=1/2$, $\tau=5\sqrt{2}$, $\sigma=2/5$, and $\kappa=1/5$. The Fourier amplitudes do not tend to zero,
indicating nonconvergence of the Fourier series and nonexistence of a (sufficiently regular) discrete travelling wave.}
\vspace{-1ex}
\end{center}
\end{figure}

Substituting the Fourier series
\begin{equation}
\label{FS}
\varphi(\xi) = \sum^{\infty}_{n = -\infty}\tilde{\varphi}_{n}e^{in\frac{\pi}{\tau}\xi}
\end{equation}
and
\begin{equation*}
h(\xi) = \sum^{\infty}_{n = -\infty}\tilde{h}_{n}e^{in\frac{\pi}{\tau}\xi}
\end{equation*}
into the difference equation \eqref{DTWE} gives
\begin{equation*}
\tilde{\varphi}_{n} = \frac{\tilde{h}_{n}}{1+d_n}
\end{equation*}
where
\begin{equation}
\label{eq:dn}
d_n = \frac{2c^{2}}{\kappa^{2}}\left(\cos\left(\frac{n\pi\kappa}{\tau}\right) -1\right) - \frac{2}{\sigma^{2}}\left(\cos\left(\frac{n\pi\sigma}{\tau}\right) -1\right)
\end{equation}
and
\begin{equation*}
\tilde{h}_{n} = \frac{1}{in{\pi}}\left((-1)^{n} - 1\right).
\end{equation*}
If the Fourier series \eqref{FS} converges, then it determines a discrete periodic travelling wave;
if not, there is no discrete periodic travelling wave of sufficient regularity to have a Fourier series.
 We make observations on two cases.
\begin{itemize}
\item[Case 1.] If $\sigma$, $\kappa$, and $\tau$ are all rationally related
(i.e. if $\sigma/\tau$, $\kappa/\tau\in\mathbb{Q}$), then the series \eqref{FS} converges, 
assuming $d_n\ne 0$ for all odd $n$. In this case \eqref{FS} is a finite sum of
Fourier series of the form $\sum_{n=-\infty}^\infty e^{i n \pi \xi/\tau}/(a + b n)$ ($a,b\in\mathbb{Z}$)
which is the Fourier series of a piecewise-constant function. Some examples in which $\varphi(\xi)$
is approximated simply by truncating the sum are shown in Figure \ref{fig:rrat}.
The expected Gibbs phenomenon is seen, together with nonmonotonic behaviour of the solutions,
similar to that seen in the heteroclinic case considered previously.
\item[Case 2.] If either or both of $\sigma/\tau$ or $\kappa/\tau$ is irrational,
then the sum \eqref{FS} suffers from small denominators. A simple statistical model,
in which $e^{i n \pi \kappa/\tau}$ are uniformly distributed on the unit circle
independently of $e^{i n \pi \sigma/\tau}$, predicts that $(1-d_n)^{-1} = \mathcal{O}(n)$
and thus $\tilde\varphi_n = \mathcal{O}(1)$ and the series \eqref{FS} does not converge.
This prediction is borne out in examples. For example, Figure \ref{fig:PhiNonConverge} shows
$|\tilde\varphi_n|$ for $c=1/2$, $\kappa=1/5$, $\sigma=2/5$, and $\tau = 5\sqrt{2}$.
While for moderate $n$ the Fourier coefficients appear to be decreasing in magnitude,
for $n\gtrsim10^7$ the asymptotic behaviour $\tilde\varphi_n = \mathcal{O}(1)$ is clear.
\end{itemize}

We conclude that the McKean caricature {\em does} allow the explicit calculation of travelling
waves for discretizations of Hamiltonian PDEs.
The nonpreservation of travelling waves when $\sigma$, $\kappa$, and $\tau$ are not rationally
related is not too suprising in view of the known sensitivity of periodic orbits of Hamiltonian ODEs
to nonsmooth perturbations. Indeed, it appears that the McKean nonlinearity, despite being discontinuous
at $u=a$, only just fails to be smooth enough---an extra factor of $\frac{1}{n}$ in $\tilde h_n$,
which would arise if the forcing were merely continuous, would be enough to make the series
\eqref{FS} converge for all parameter values.
This motivates our study in the next section of a continuous, piecewise linear nonlinearity.

\section{Continuous but non-smooth $V^{\prime}$: Sawtooth nonlinearity} \label{sec:sawtooth}

\begin{figure}[h]
\centering
\includegraphics[width=0.4\textwidth]{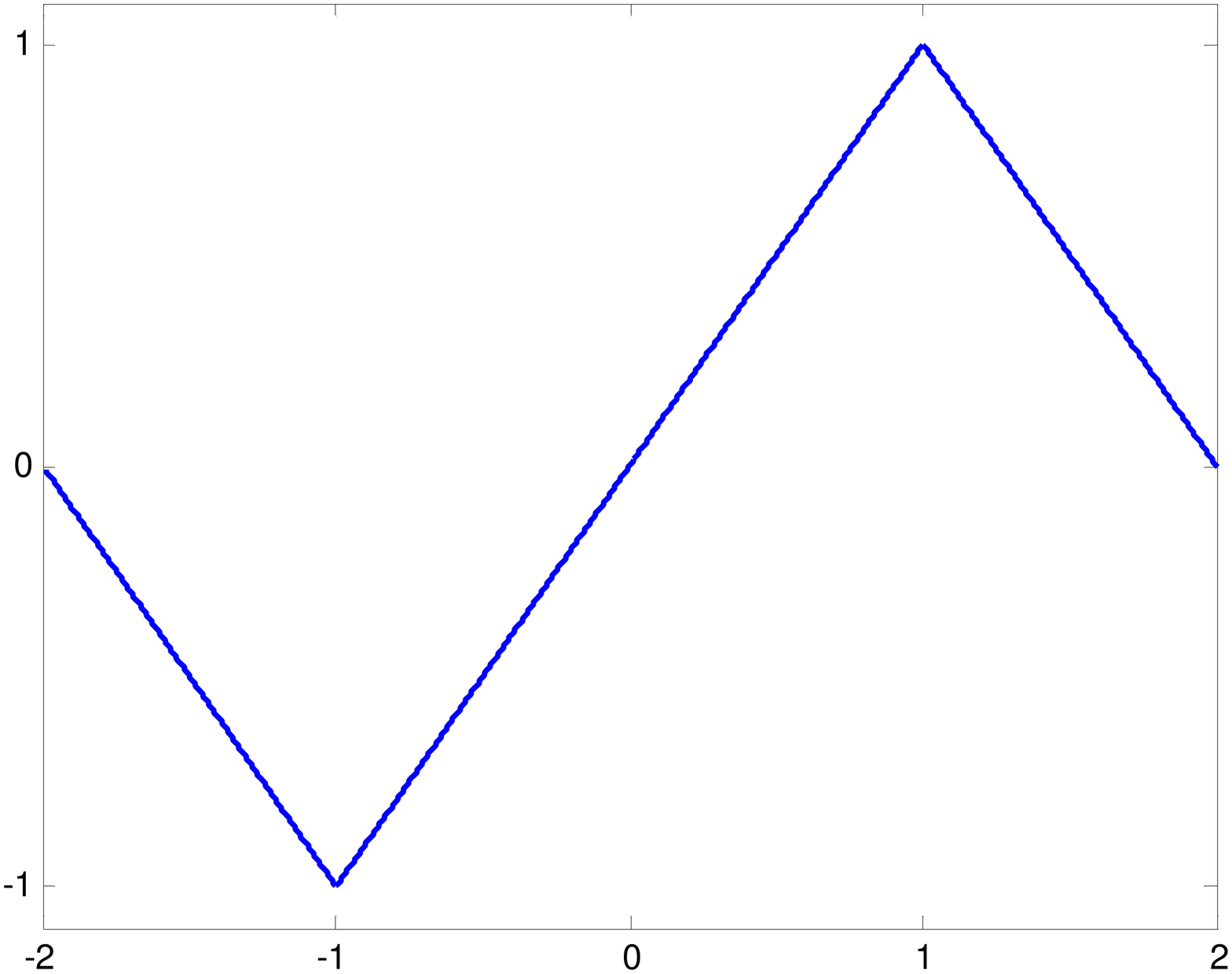}
\caption{Sawtooth function $-V^{\prime}(\varphi)$. \label{NewNon}}
\end{figure}

\noindent
We now consider the continuous piecewise linear approximation of the cubic
with domain $\varphi\in(-2,2)$
given by
\begin{equation}
\label{Sawtooth}
V^{\prime}(\varphi) = \left\{
\begin{array}{rl}
\varphi + 2, & -2<\varphi<-1\\
-\varphi, & -1<\varphi<1\\
\varphi - 2, & 1<\varphi<2
\end{array} \right.
\end{equation}
which is shown in Figure~\ref{NewNon}.
The phase portrait of the ODE \eqref{eq:reduced}
with nonlinearity \eqref{Sawtooth} is given in Figure~\ref{SawtoothPhaseAndExpected},
and shows that the equation has heteroclinic and periodic travelling waves.
We choose to concentrate on periodic travelling waves.

\subsection{Continuous travelling waves}

\begin{figure}[h]
\centering
{\label{SawtoothPhaseLabels}\includegraphics[width=0.49\textwidth]{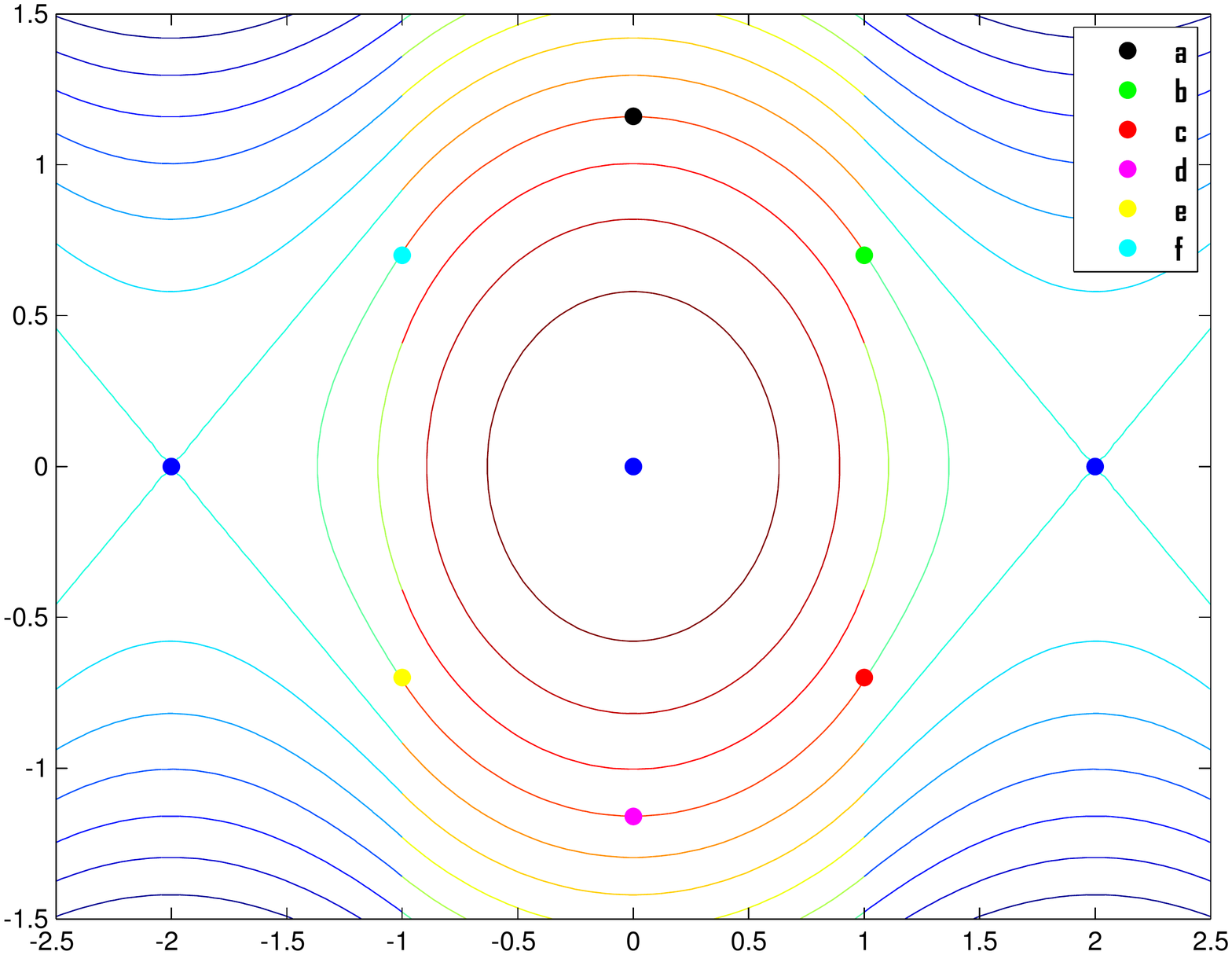}}
{\label{SawtoothExpectedSolution}\includegraphics[width=0.49\textwidth]{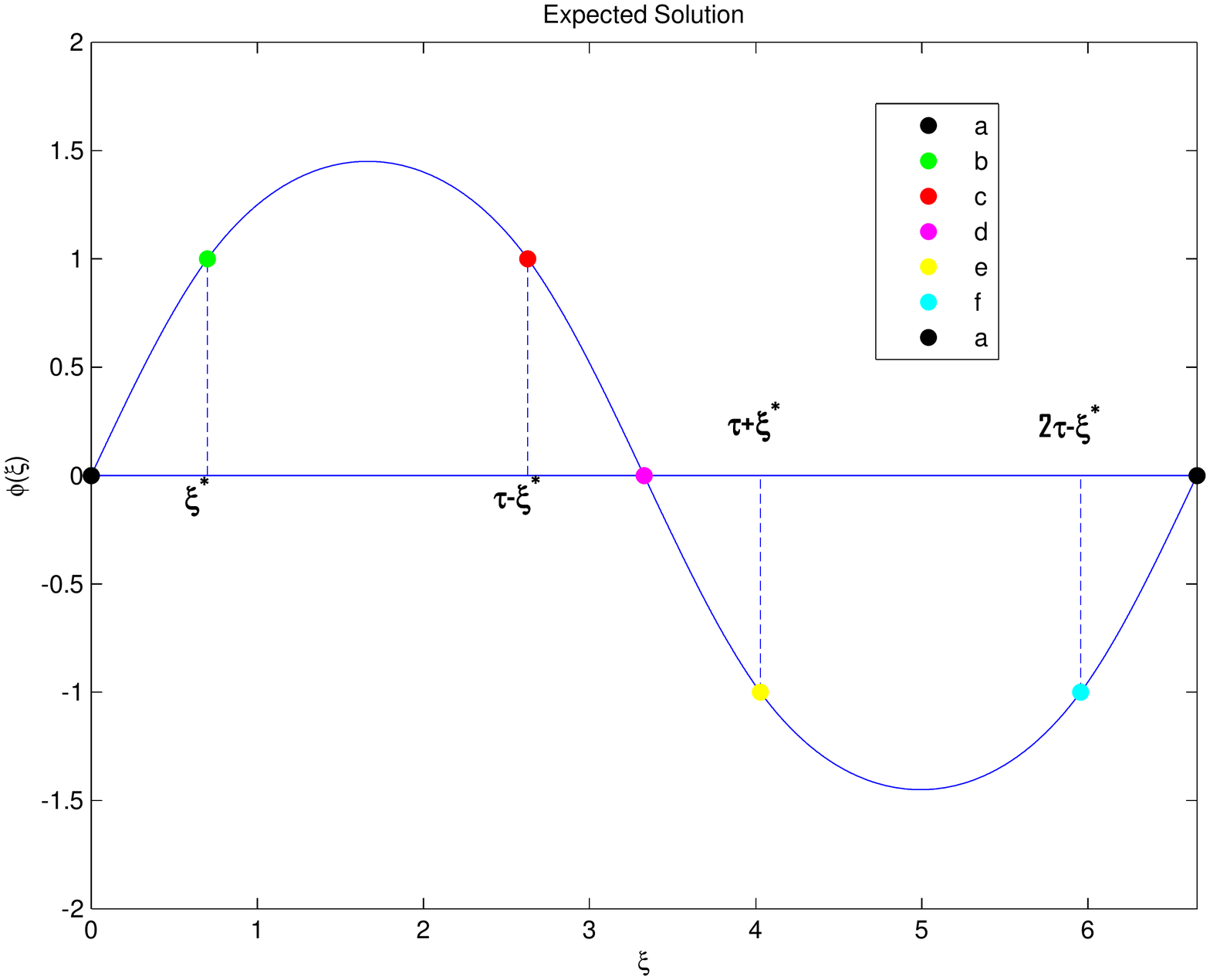}}
\caption{Phase portrait (left) and expected solution (right) of the nonlinear wave equation with sawtooth nonlinearity.}
\label{SawtoothPhaseAndExpected}
\end{figure}

Periodic solutions with amplitude greater than 1 will have the structure illustrated
in Figure \ref{SawtoothPhaseAndExpected}. We fix the period to be $2\tau$ and introduce
the 6 marked points a,b,\dots,f. By translation symmetry, we are free to take
$\varphi(0)=0$ (point `a'). We define $\xi^*$ as the least value of $\xi$ such that
$\varphi(\xi^*)=1$ (point `b'). By symmetry, we then have $\varphi(\tau-\xi^*)=1$
(point `c'), $\varphi(\tau)=0$ (point `d'), $\varphi(\tau+\xi^*)=-1$ (point `e'), 
$\varphi(2\tau-\xi^*)=-1$ (point `f'), and $\varphi(2\tau)=0$ (point `a').

Therefore the nonlinearity becomes
\begin{equation}
\label{Sawtooth2}
f(\varphi) := -V'(\varphi)=\left\{
\begin{array}{rl}
\varphi, & 0<\xi<\xi^* \\
-\varphi+2, & \xi^*<\xi<\tau-\xi^*\\
\varphi, & \tau-\xi^*<\xi<\tau+\xi^* \\
-\varphi-2, & \tau+\xi^*<\xi<2\tau-\xi^* \\
\varphi, & 2\tau-\xi^*<\xi<2\tau.
\end{array} \right.
\end{equation}
On the intervals $(0,\xi^*)$, $(\tau-\xi^*,\tau+\xi^*)$, and $(2\tau-\xi^*,2\tau)$, the differential equation is
\begin{equation*}
(c^2-1)\varphi''(\xi)=\varphi(\xi)
\end{equation*}
with  general solution
\begin{equation*}
\varphi(\xi)=c_1\cos\left(\frac{\xi}{\sqrt{1-c^2}}\right)+c_2\sin\left(\frac{\xi}{\sqrt{1-c^2}}\right)
\end{equation*}
where $c_1,c_2$ are different constants for each interval determined by the interval's boundary conditions
\begin{align*}
0<\xi<\xi^*, \text{  } & \varphi(0)=0, \text{  } \varphi(\xi^*)=1, \\
\tau-\xi^*<\xi<\tau+\xi^*, \text{  } & \varphi(\tau-\xi^*)=1, \text{  } \varphi(\tau+\xi^*)=-1, \\
2\tau-\xi^*<\xi<2\tau, \text{  } & \varphi(2\tau-\xi^*)=-1, \text{  } \varphi(2\tau)=0.
\end{align*}
On  $(\xi^*,\tau-\xi^*)$ we have
\begin{equation*}
(c^2-1)\varphi''(\xi)=-\varphi(\xi)+2
\end{equation*}
with  general solution
\begin{equation*}
\varphi(\xi)=c_3 e^{\frac{\xi}{\sqrt{1-c^2}}}+c_4 e^{\frac{-\xi}{\sqrt{1-c^2}}} + 2
\end{equation*}
where $c_3,c_4$ are constants determined by the boundary conditions $\varphi(\xi^*)=\varphi(\tau-\xi^*)=1$.
On $(\tau+\xi^*,2\tau-\xi^*)$ we have the general solution
\begin{equation*}
\varphi(\xi)=c_5 e^{\frac{\xi}{\sqrt{1-c^2}}}+c_6 e^{\frac{-\xi}{\sqrt{1-c^2}}} - 2
\end{equation*}
where $c_5,c_6$ are constants determined by the boundary conditions $\varphi(\tau+\xi^*)=\varphi(2\tau-\xi^*)=-1$.
Figure 
\ref{SawtoothIrrat} shows solutions to the nonlinear wave equation
with sawtooth nonlinearity for different values of $\xi^*$ and fixed $c$.



\subsection{Discrete travelling waves}
\label{DiscreteSolution}

We now examine the periodic  solutions of the discrete travelling wave equation \eqref{DTWE}
with sawtooth nonlinearity \eqref{Sawtooth}.
As in the continuous case, for solutions of the form shown in Figure \ref{SawtoothPhaseAndExpected},
the equation becomes linear in $\varphi$
where $-V^{\prime}(\varphi)$ is given in \eqref{Sawtooth2}.

\begin{figure}[h]
\centering
\includegraphics[width=0.9\textwidth]{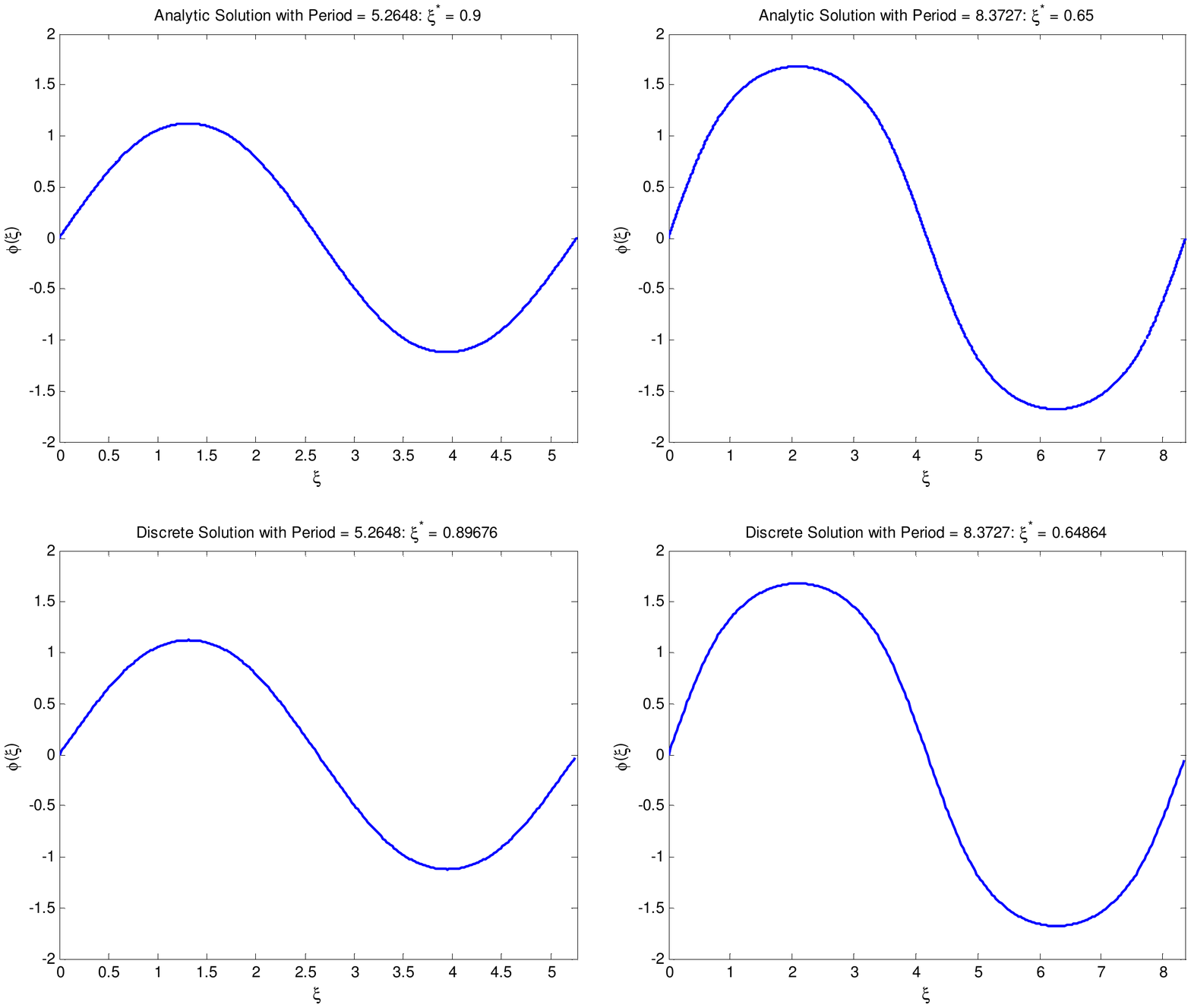}
\caption{Comparison of the analytic solution and the discrete solution of the
nonlinear wave equation with sawtooth nonlinearity, for fixed periods and fixed $c$,
for irrational values of $\frac{\sigma}{\kappa}$. \label{SawtoothIrrat}}
\end{figure}

Following Section~\ref{PeriodicSolutions}, we expand both sides of  \eqref{DTWE}  in a Fourier series, obtaining
\begin{equation}
\label{SawtoothFS}
\tilde{\varphi}_{n}d_n  = \tilde{f}_{n}.
\end{equation}
where $d_n$ is given in \eqref{eq:dn}.
Since $f(\varphi(\xi))$ is a periodic function with fixed period, $2\tau$,
for fixed $c$ and each of the five intervals in the nonlinearity \eqref{Sawtooth2},
the Fourier coefficients $\tilde{f}_{n}$ can be found explicitly and are given by
\begin{align*}
\tilde{f}_{n}&=\frac{1}{2\tau}\int^{2\tau}_{0}f(\varphi(\xi))e^{-in\frac{\pi}{\tau}\xi}d\xi\\
&=\frac{1}{2\tau}\left[\int^{\xi^*}_{0}\varphi(\xi)e^{-in\frac{\pi}{\tau}\xi}d\xi+\int^{\tau-\xi^*}_{\xi^*}(-\varphi(\xi)+2)e^{-in\frac{\pi}{\tau}\xi}d\xi+ \right.  \\
&\hspace{-5pt}\left.\int^{\tau+\xi^*}_{\tau-\xi^*}\varphi(\xi)e^{-in\frac{\pi}{\tau}\xi}d\xi+\int^{2\tau-\xi^*}_{\tau+\xi^*}(-\varphi(\xi)-2)e^{-in\frac{\pi}{\tau}\xi}d\xi+\int^{2\tau}_{2\tau-\xi^*}(\varphi(\xi))e^{-in\frac{\pi}{\tau}\xi}d\xi\right]\\
&=\sum^{\infty}_{k=-\infty}A_{kn}\tilde{\varphi}_{n}+b_n
\end{align*}
where
\begin{align*}
A_{kn} &= \frac{1+(-1)^{(k-n)}}
{i(k-n)\pi}\left(e^{i(k-n)\frac{\pi}{\tau}\xi^*}-e^{-i(k-n)\frac{\pi}{\tau}\xi^*}\right)\\
b_n&=\frac{1-(-1)^{n}}{in\pi}\left(e^{-in\frac{\pi}{\tau}\xi^*}+e^{in\frac{\pi}{\tau}\xi^*}\right).
\end{align*}
Substituting this into \eqref{SawtoothFS} and solving for $\tilde{\varphi}_{n}$ we get,
\begin{equation}
\label{SawtoothFSC}
\tilde{\varphi}_{n}=(\mathrm{diag}(d)-A)^{-1}b.
\end{equation}

Thus far it appears that finding a travelling wave only requires solving a linear system.
However, there is one final condition that needs to be imposed, and this one is nonlinear.
This is the compatibility condition $\varphi(\xi^*)=1$, which has to be solved numerically.
The piecewise-linear force has reduced the problem to
solving a scalar nonlinear equation in one variable. Our algorithm for finding discrete travelling waves is thus:
(i) Fix the parameters $c$, $\sigma$, $\kappa$, and $\tau$.
(ii) Choose a suitable initial number $N$ of Fourier modes with which to represent $\varphi$,
and an initial guess for $\xi^*$ equal to that corresponding to the analytic solution with the same value of $\tau$.
(iii) Solve the equation $\varphi(\xi^*)=1$ numerically. This involves repeatedly forming and solving an $N\times N$
linear system. (iv) Increase $N$ and examine the convergence of the Fourier series.
In practice, we found that the nonlinear solve was extremely easy.
The main drawback of the algorithm is that the linear solves limit the value of $N$ that can be used.

Discrete numerical solutions are plotted and compared with the analytic solution in Figure 
\ref{SawtoothIrrat}.
It appears that the discrete numerical solutions are very similar to the analytic solution.
This was confirmed by plotting the discrete solutions on top of the analytic solutions for fixed periods,
varying the ratio $\sigma/\kappa$. In all cases tried, for both rational and irrational $\sigma/\kappa$, 
the discrete numerical solution seemed to match the
corresponding analytic solution closely. No wiggles or lack of convergence were observed.
Thus, the piecewise linear model appears to be a promising candidate for further study.

\section{Smooth $V^{\prime}$: Sine nonlinearity} \label{sec:smooth}

Our final choice of nonlinearity is a smooth nonlinear function. We choose the sine function,
so that our nonlinear wave equation becomes the sine-Gordon equation
\begin{equation}
\label{sineGordon}
u_{tt}=u_{xx}-\sin{u}.
\end{equation}
Thus, we take $V^{\prime}(\varphi) = \sin(\varphi)$ in equations \eqref{eq:reduced} and \eqref{DTWE}.
Again, there are periodic and heteroclinic travelling waves.

We use a pseudospectral Newton continuation method to find periodic travelling waves.
We will be working in both real space and Fourier space so make the following definitions.
Let $-\sin(\varphi(\xi))=f(\varphi(\xi))$. Discrete periodic travelling waves will be approximated by the Fourier series
\begin{equation*}
\varphi(\xi)=\sum_{n=-N}^{N}\tilde{\varphi}_{n}e^{in\frac{\pi}{\tau}\xi}
\end{equation*}
where $T = 2\tau$ is the period of the  solution.
As in section~\ref{DiscreteSolution} we have
\begin{equation}
\label{FC}
D\tilde{\varphi}=\tilde{f}
\end{equation}
where $D=\mathrm{diag}(d)$, with $d$ defined in \eqref{eq:dn}.
In real space,
\begin{equation}
\label{SmoothNewton}
L\varphi=f(\varphi)
\end{equation}
where $L=\mathcal{F}^{-1}{D}\mathcal{F}$ is a linear operator, $\varphi = (\varphi(i n \pi \xi/\tau))_{i=0}^N$,  
$\mathcal{F}\varphi=\tilde{\varphi}$, and $\mathcal{F}$ is the discrete Fourier transform.

We solve the nonlinear equation \eqref{SmoothNewton} using Newton's method,
\begin{equation*}
\varphi^{j+1}=\varphi^{j}-(L-\text{diag}(f'(\varphi^{j}))^{-1}(L\varphi^{j}-f(\varphi^{j})).
\end{equation*}
Generally we hold $\sigma$, $\kappa$, and $c$ fixed and apply continuation in the parameter $T$ starting from its minimum value $T=2\pi$.
If we have quadratic convergence of Newton's method, and the Fourier coefficients $\tilde\varphi_n$
are decreasing rapidly, indicating convergence of the Fourier series, then this numerical solution will
be  close to an exact solution of the discrete travelling wave equation \eqref{DTWE}.
This is what we observe in almost all cases, the exceptions to be noted later.

\begin{figure}
\centering
\subfloat[Orbit of leapfrog applied to \eqref{eq:reduced}]{\includegraphics[width=0.5\textwidth]{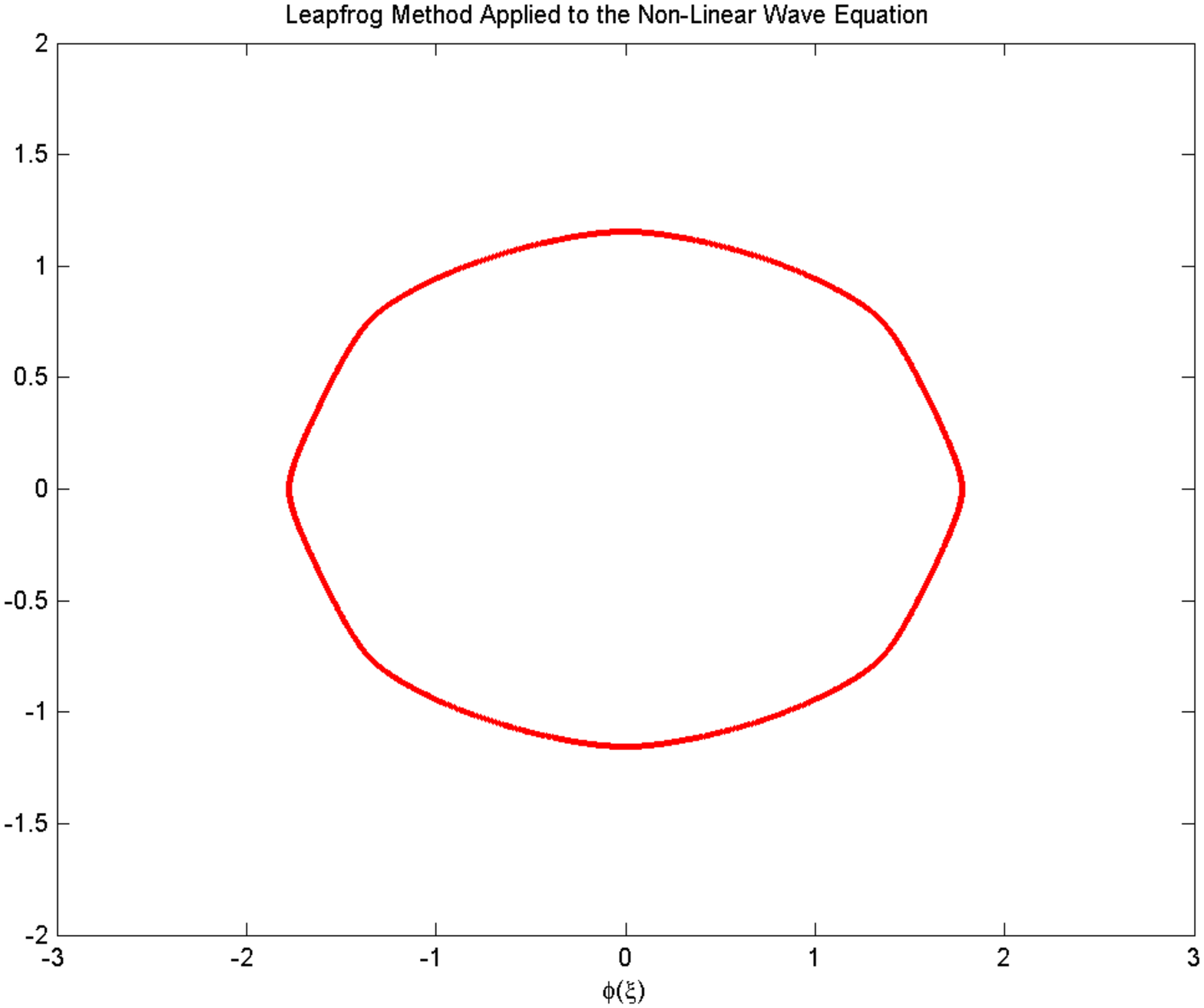}}
~
\subfloat[Red: Solution of leapfrog applied to \eqref{eq:reduced}. Blue: Numerical solution of \eqref{DTWE}]{\includegraphics[width=0.5\textwidth]{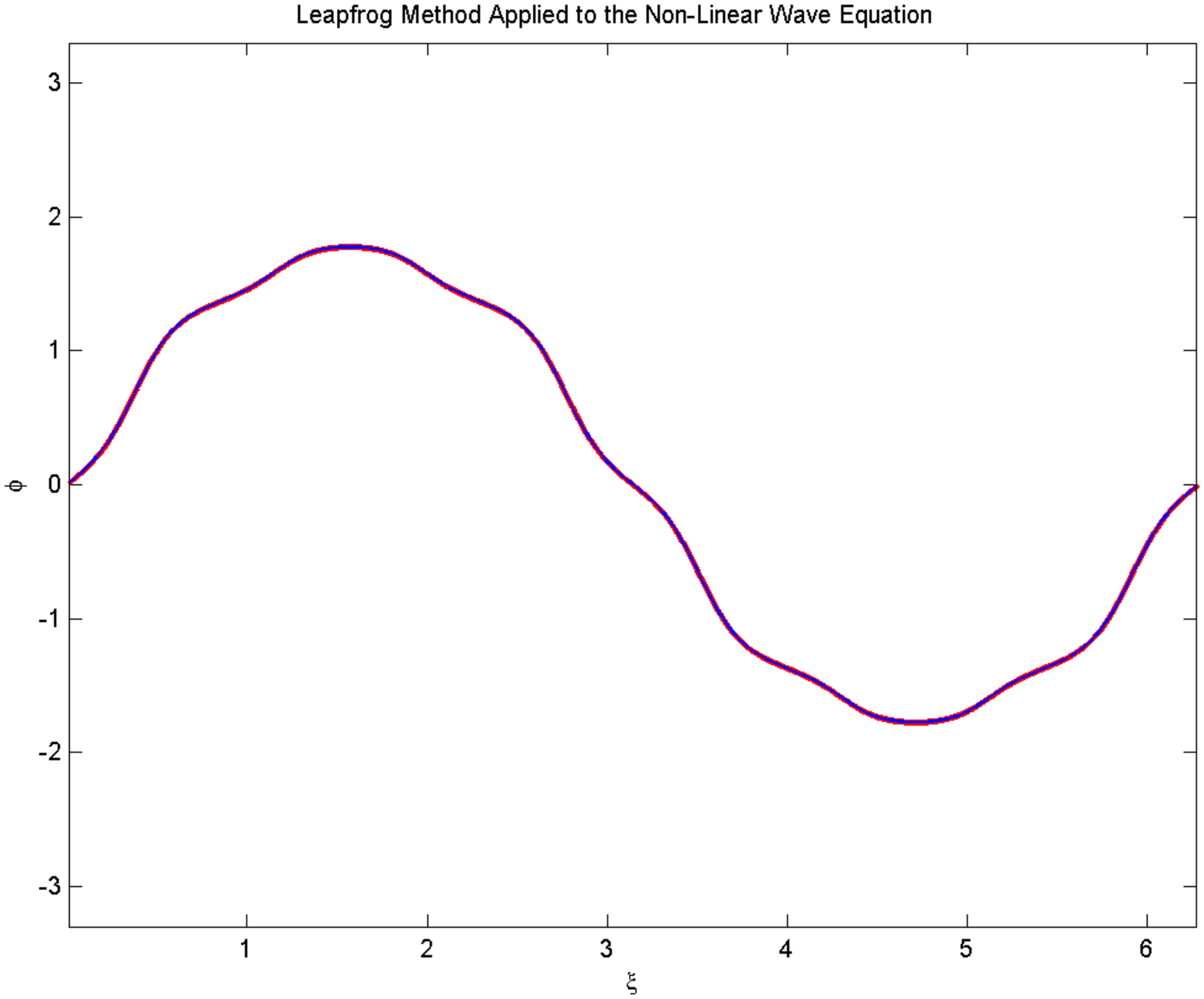}}
\caption{A comparison of solutions of \eqref{eq:reduced} and \eqref{DTWE} for $\sigma=\kappa=0.8$.}
\label{LeapfrogComp2}
\end{figure}

We first present a second, independent check that this numerical method yields
correct solutions to the discrete travelling wave equation \eqref{DTWE}.
When $\sigma=\kappa$ the discrete travelling wave equation \eqref{DTWE} is
equivalent to finding periodic solutions of the leapfrog method applied to the
pendulum problem, which is relatively well understood.
We compare the orbit of the leapfrog method applied to \eqref{eq:reduced} with
initial conditions taken from a solution of the pseudospectral--Newton method.
An example is shown in Figure~\ref{LeapfrogComp2}
for $\sigma=\kappa=0.8$ (known to be resonant), and the two solutions agree.

\begin{figure}
\centering
\subfloat[$T=2\pi=6.2832$]{\includegraphics[width=0.40\textwidth]{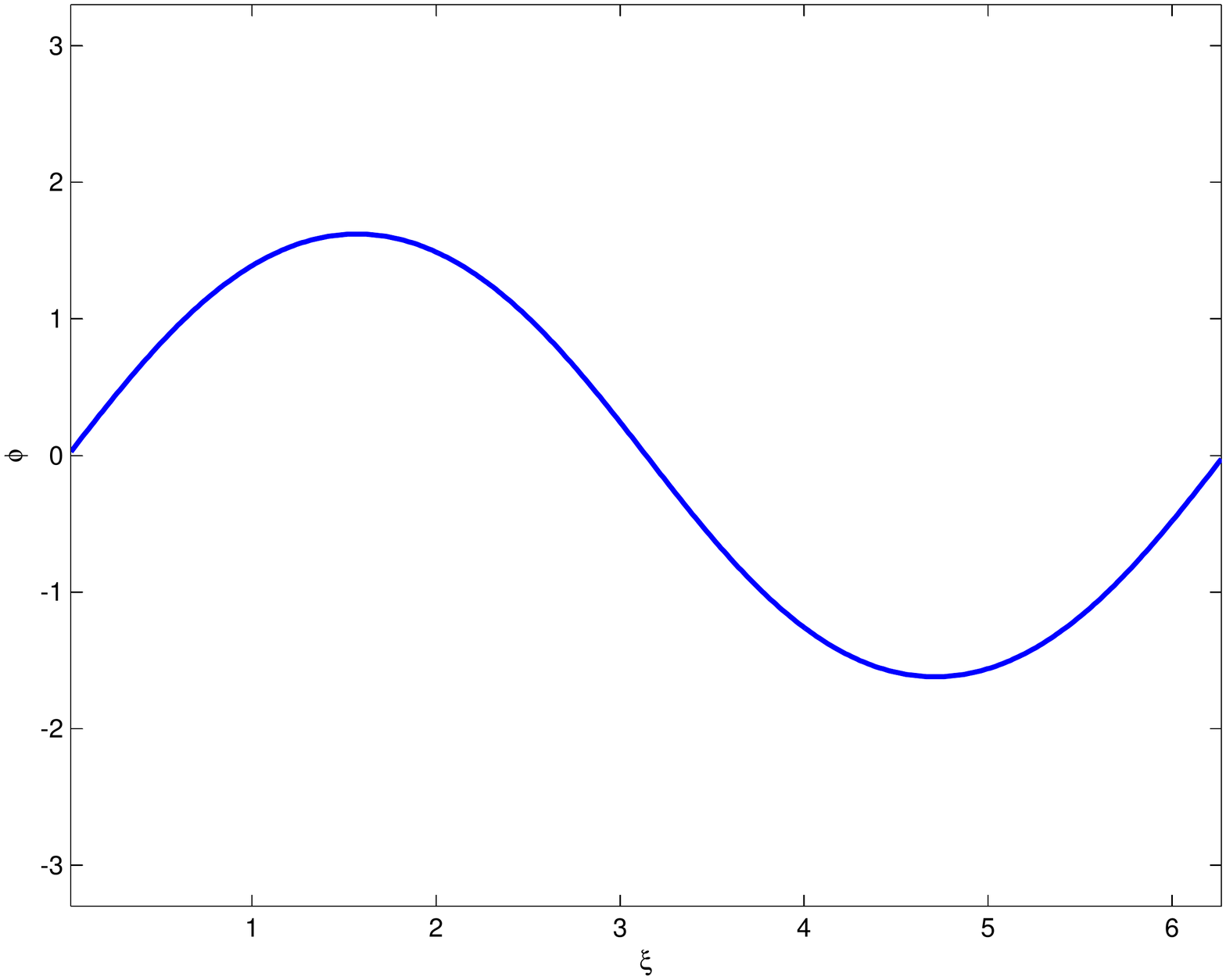}\label{ContinuationPi2}}
~
\subfloat[$T=2(\pi+4.4)=15.0832$]{\includegraphics[width=0.40\textwidth]{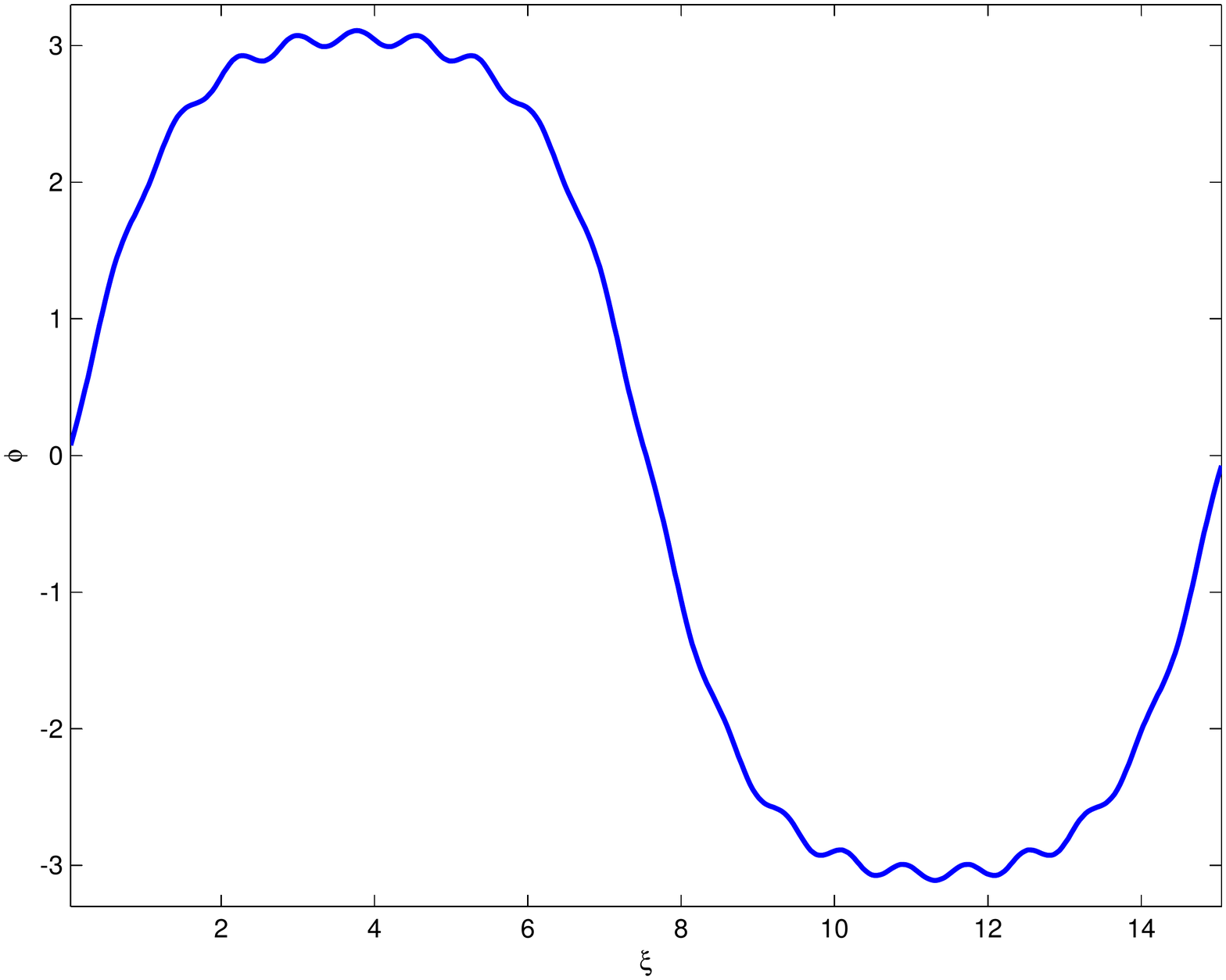}\label{Continuation76}}

\subfloat[$T=2(\pi+6.7)=19.6832$]{\includegraphics[width=0.40\textwidth]{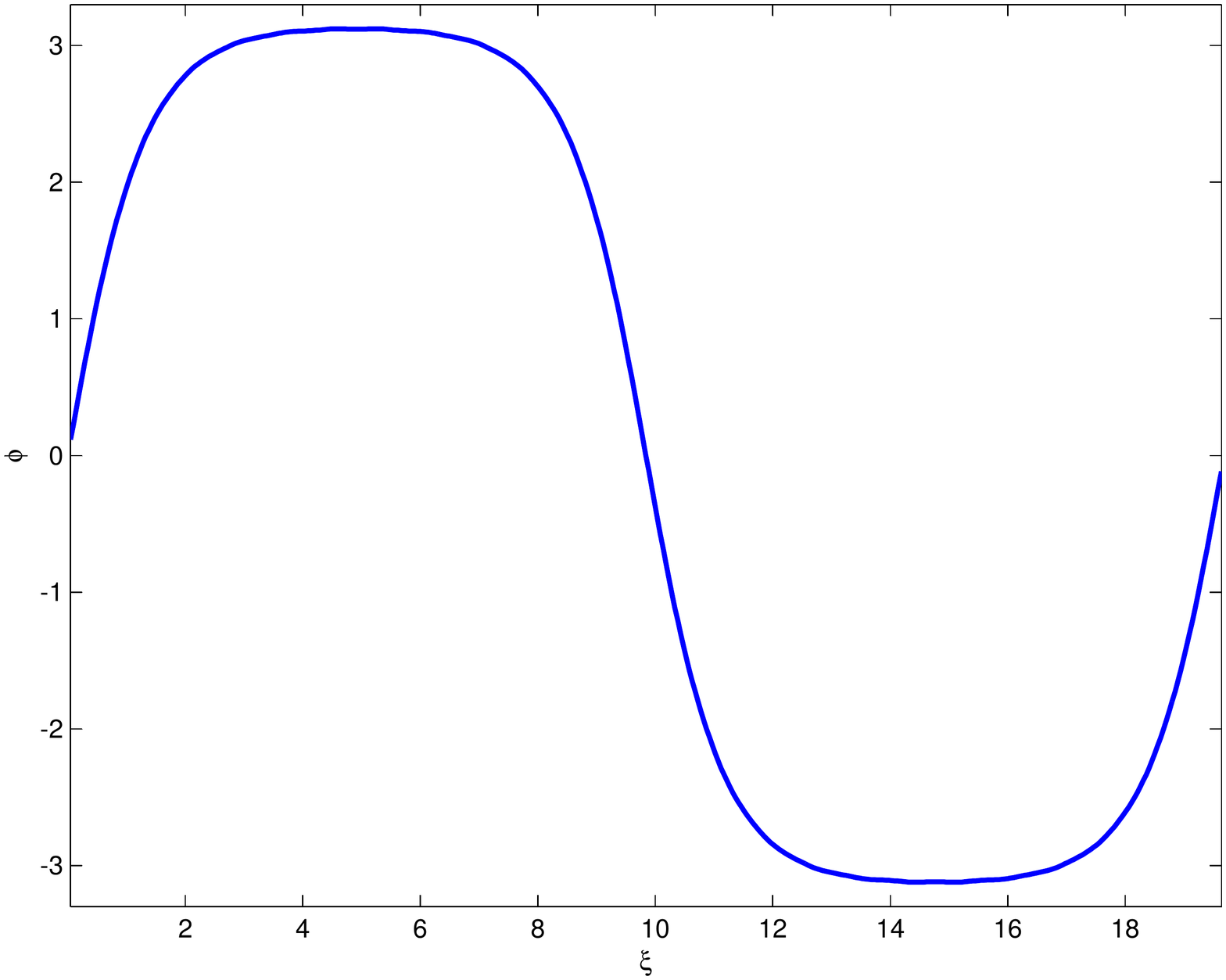}\label{Continuation99}}
~
\subfloat[$T=2(\pi+9.9)=26.0832$]{\includegraphics[width=0.40\textwidth]{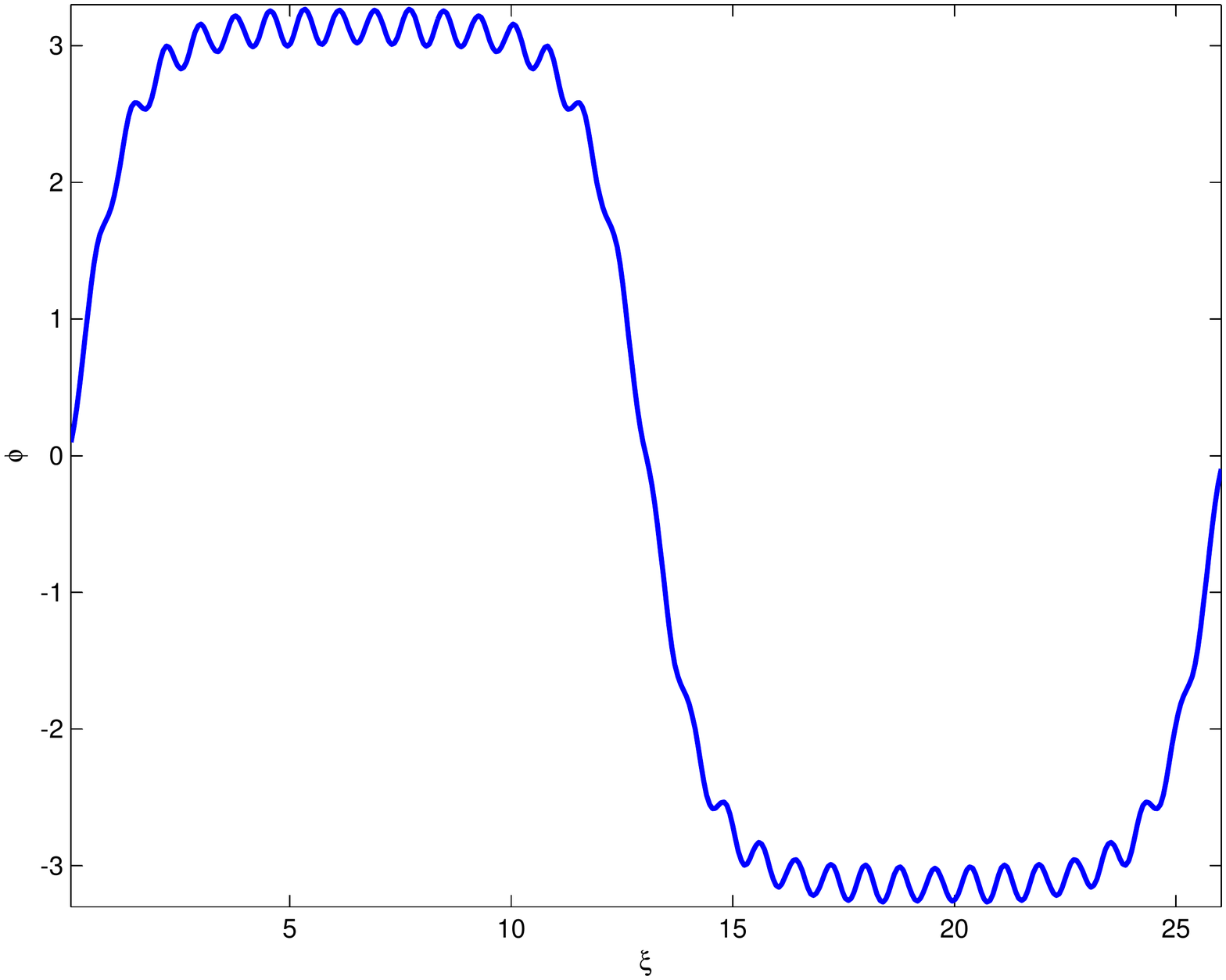}\label{Continuation131}}

\subfloat[$T=2(\pi+13.4)=33.0832$]{\includegraphics[width=0.40\textwidth]{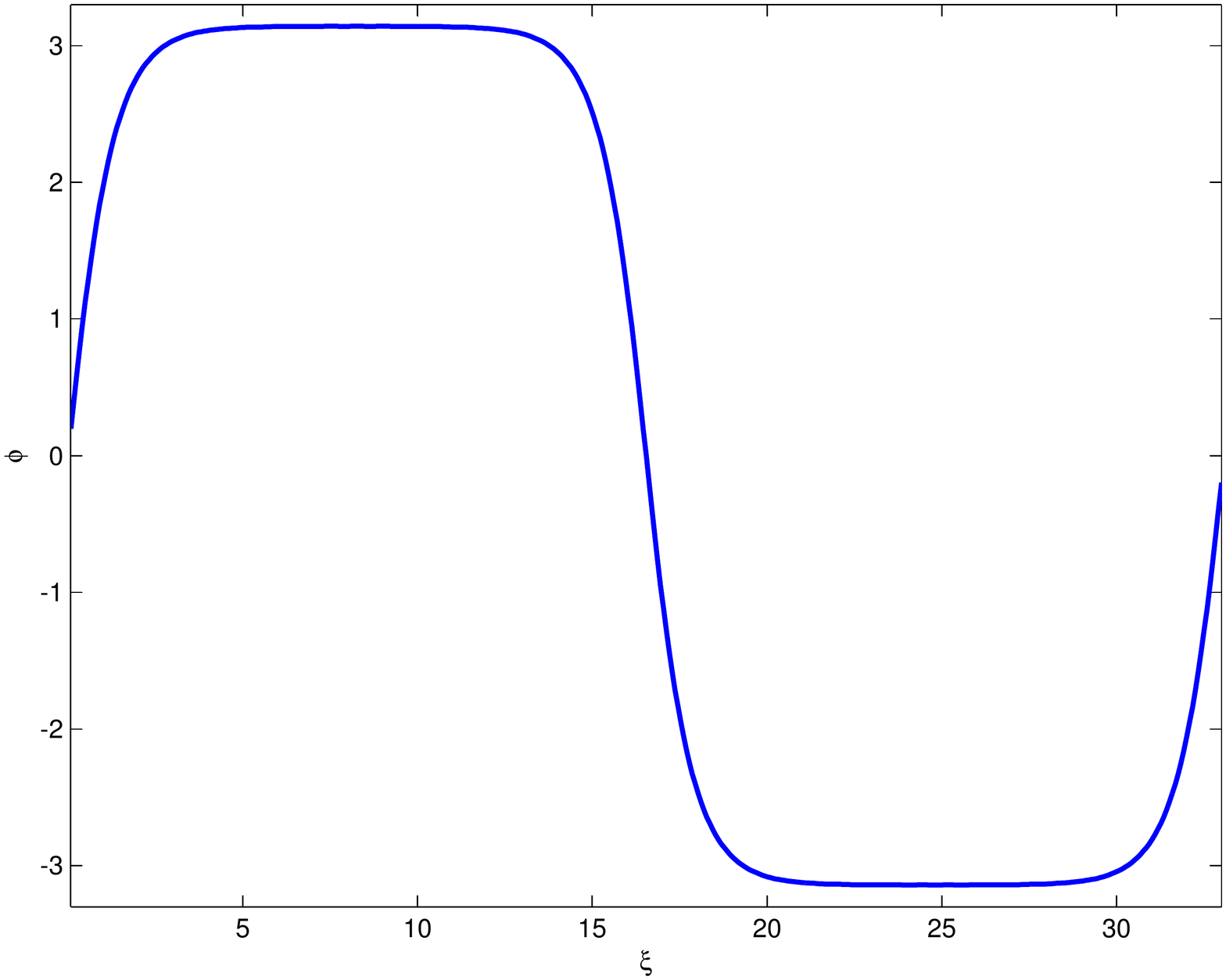}\label{Continuation166}}
~
\subfloat[$T=2(\pi+16.8)=39.8832$]{\includegraphics[width=0.40\textwidth]{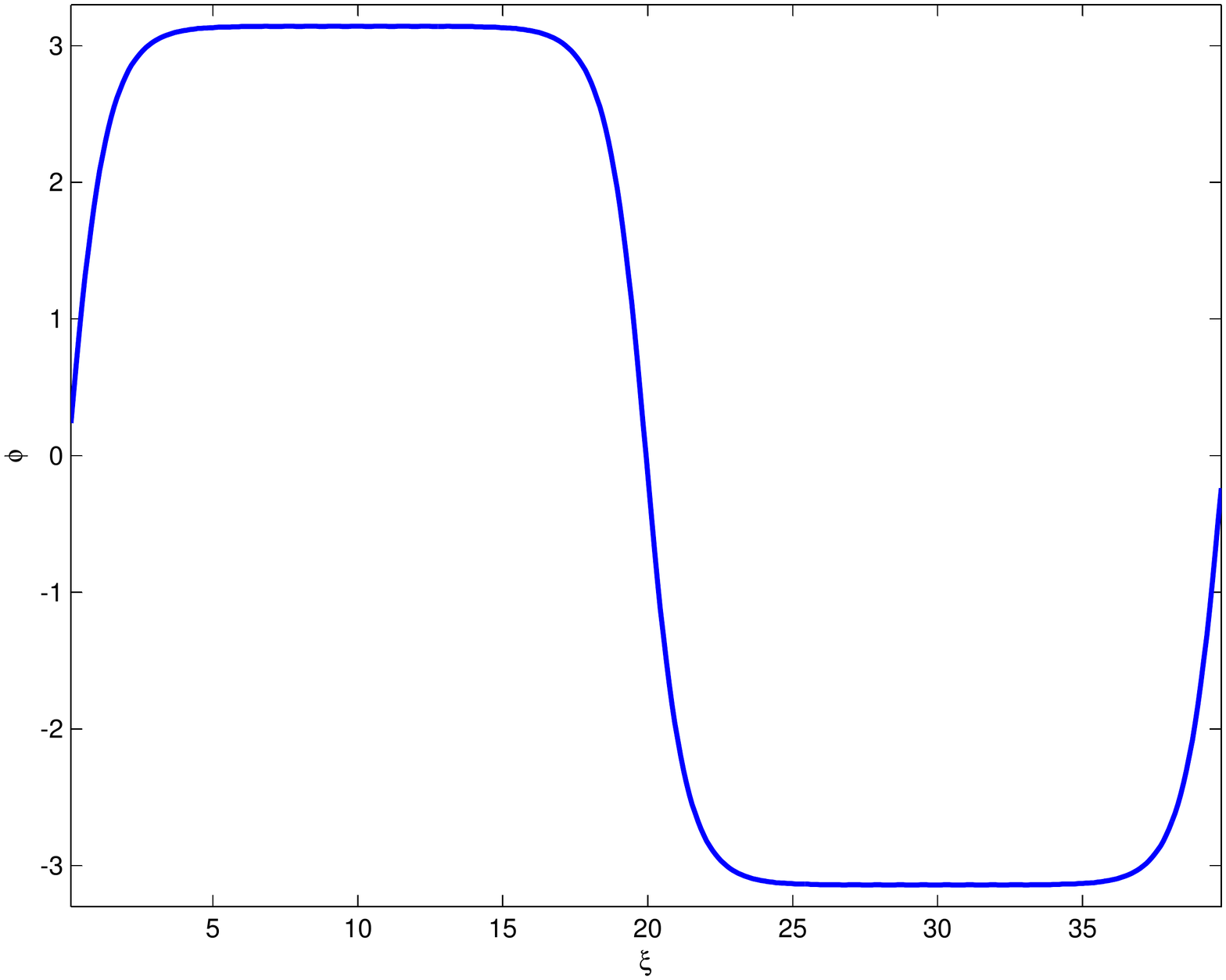}\label{Continuation20}}
\caption{Continuation in $T$ starting with $T=2\pi$ for $\sigma=1$, $\kappa=\frac{1}{\sqrt{2}}$, $c=1.3$.}
\label{Continuation2}
\end{figure}

Figure \ref{Continuation2} shows a sequence of solutions from the continuation.
There are two distinct types of solution: those that are smooth we call non-resonant,
those that possess noticeable wiggles we call resonant.

\begin{figure}[h]
\centering
\subfloat[$T=2(\pi+0.49)=7.2632$]{\includegraphics[width=0.40\textwidth]{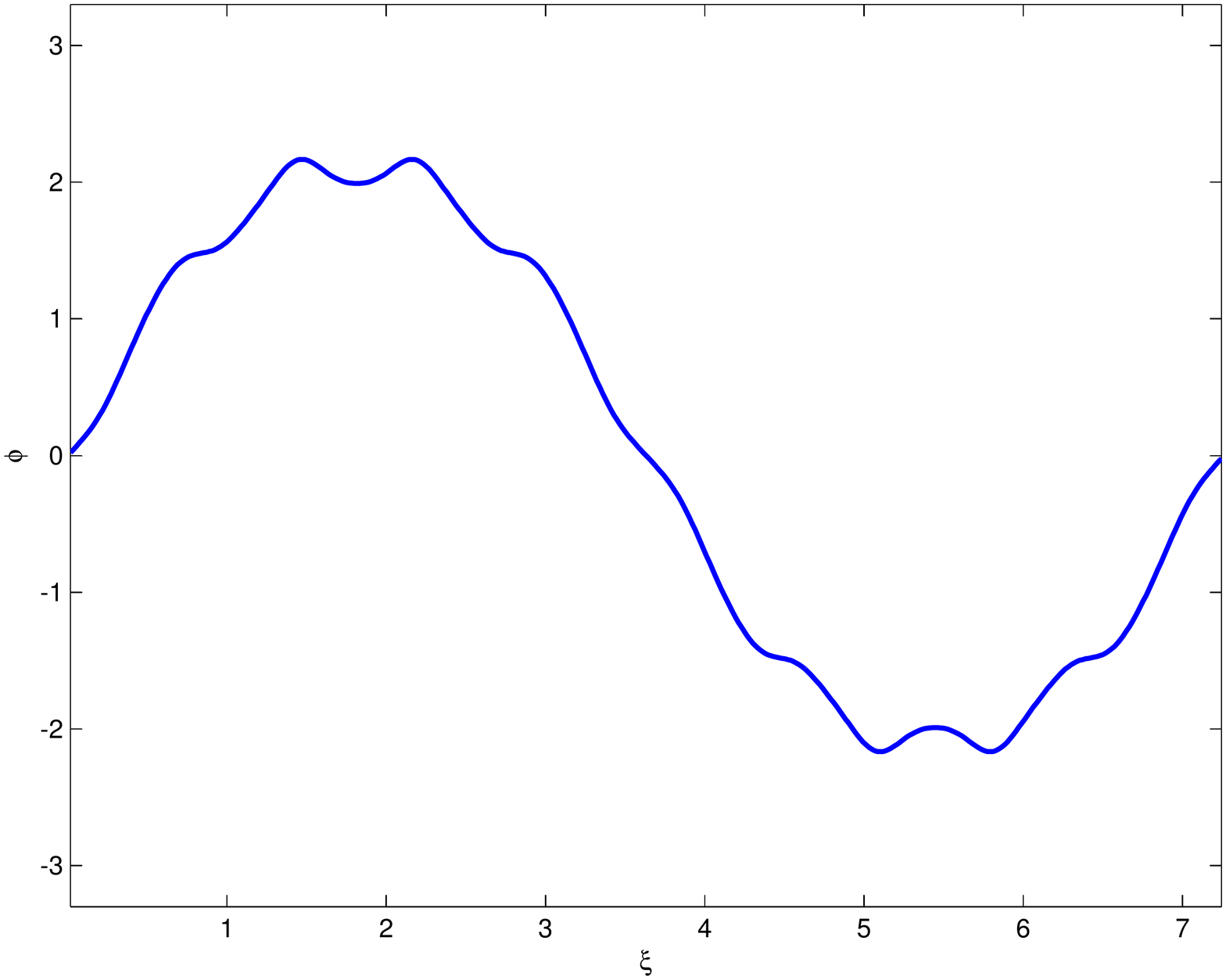}\label{Continuation3642}}
~
\subfloat[Resonances]{\includegraphics[width=0.40\textwidth]{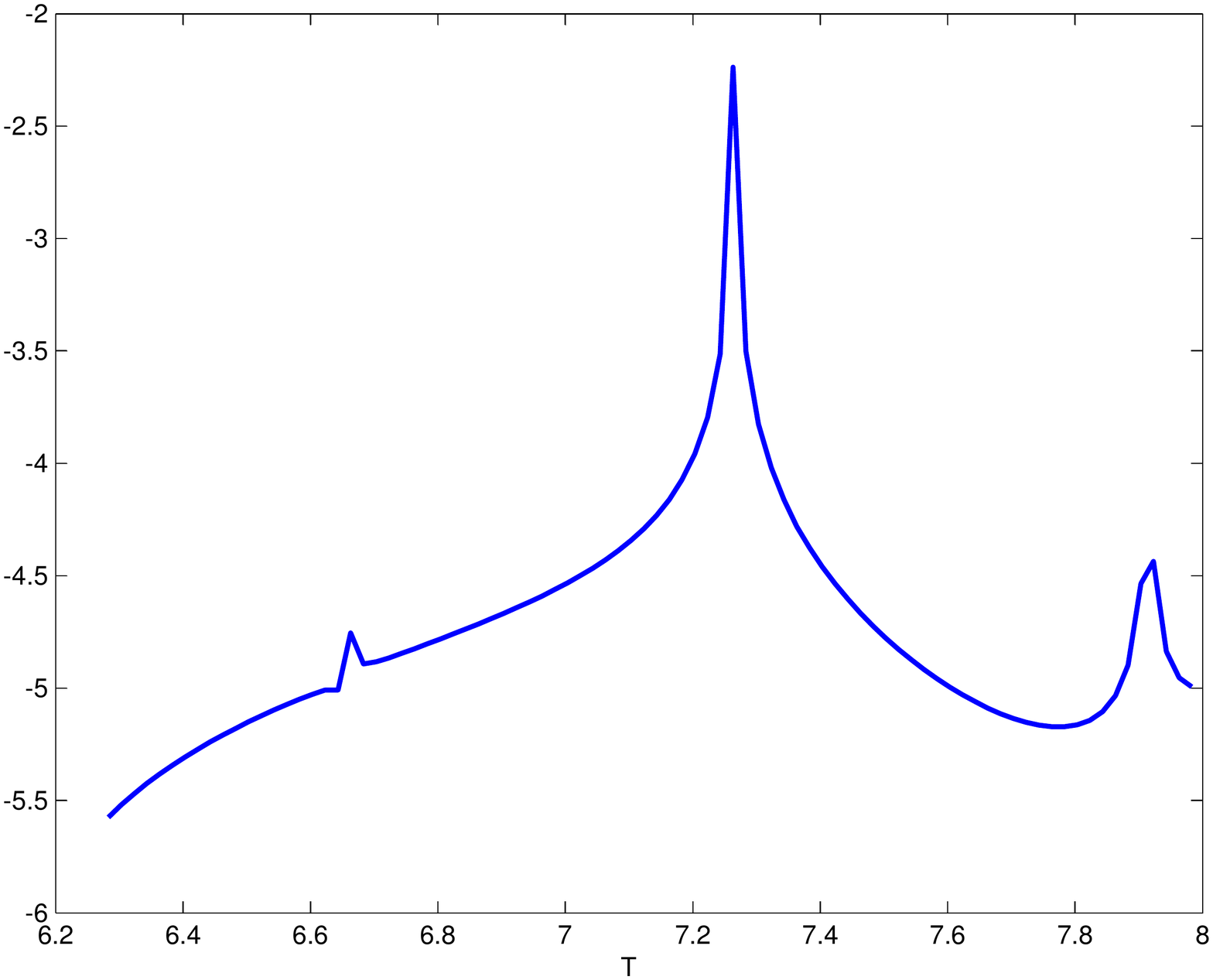}\label{Resonance4}}
\caption{A short continuation in $T$ showing one major resonance peak. The $y$-axis gives the size of the resonance defined by \eqref{ResSize}. \label{Resonance44}}
\end{figure}

To quantify the resonances we use an ad-hoc measure based on the observation that
smooth solutions for large $T$ are very flat near $\xi = T/4$. After some experimentation the function
\begin{equation}
\label{ResSize}
R :=\log_{10}(\max_{0.2T<\xi<0.3T}\left|d''\right|-\min_{0.2T<\xi<0.3T}\left|d''\right|).
\end{equation}
was found to be a useful measure of resonance amplitude and was able to detect resonances.
The resonance amplitude $R$  against the period of the solution
is shown in Figure~\ref{Resonance4}, in which the major peak showing
corresponds to the solution shown in Figure \ref{Continuation3642}.

\begin{figure}[h]
\centering
\includegraphics[width=0.49\textwidth]{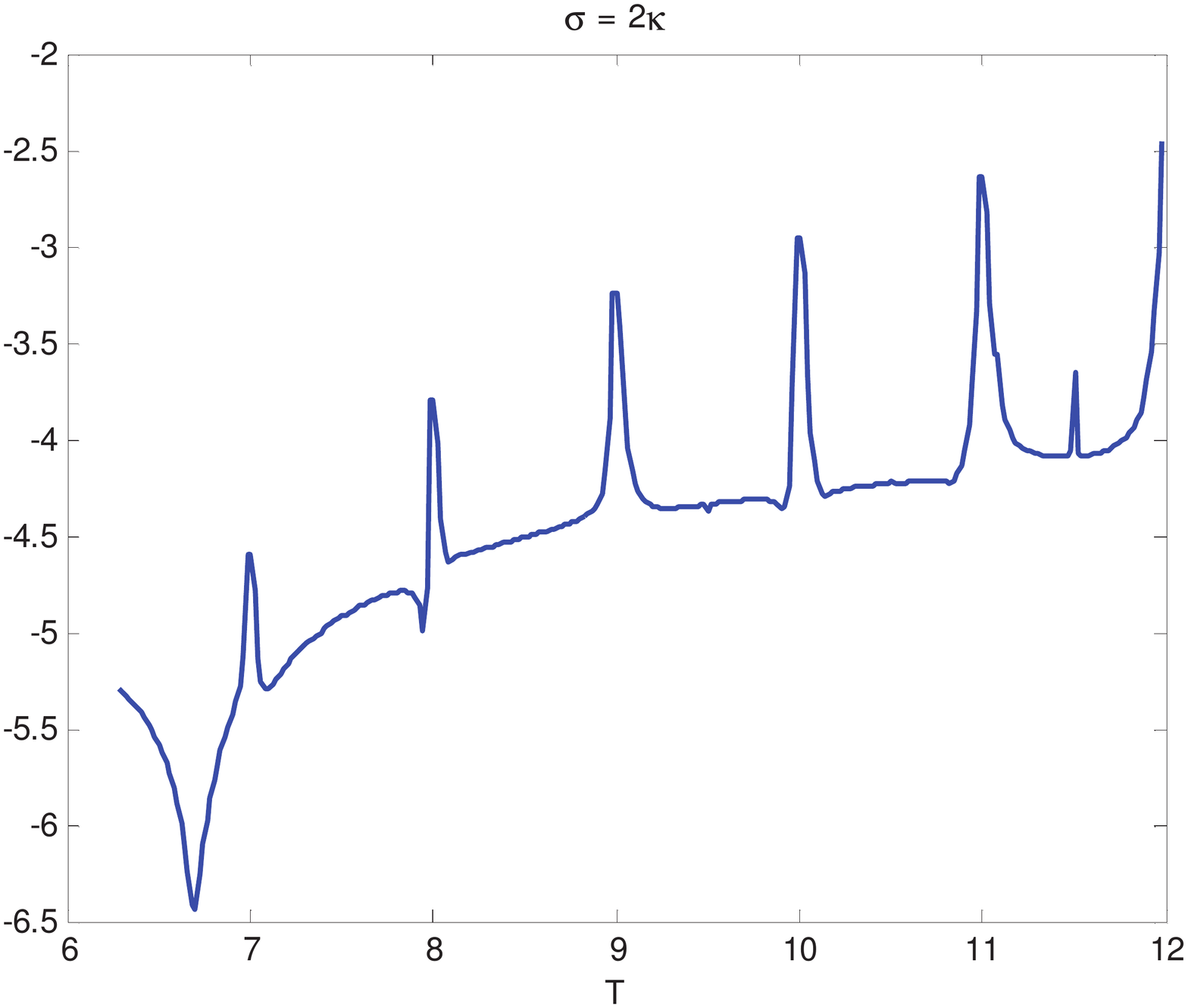}
\includegraphics[width=0.49\textwidth]{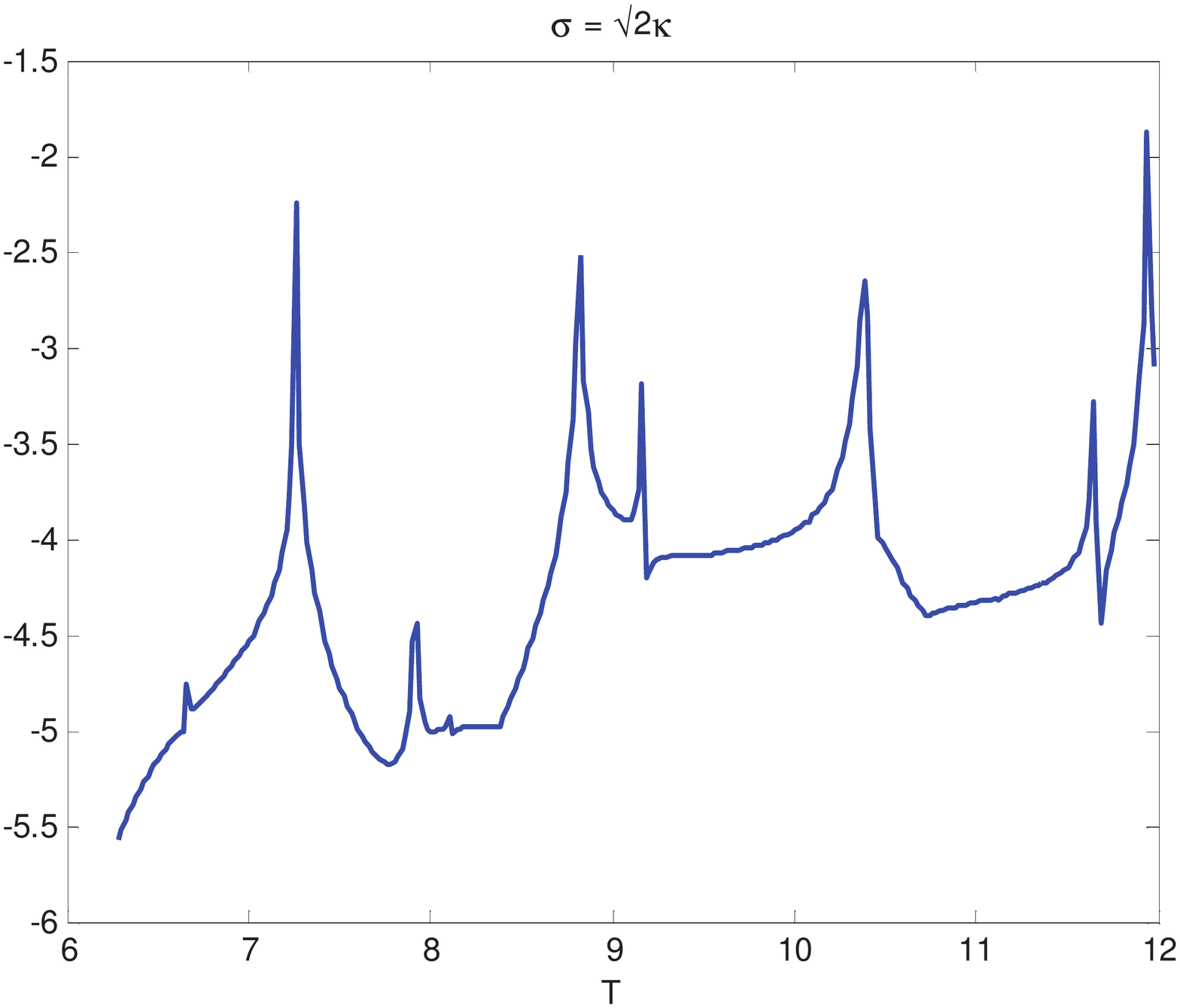}
\vspace{-10pt}
\caption{Comparison of the resonances for a rational ratio of $\frac{\sigma}{\kappa}$ and an irrational one. 
The $y$-axis gives the size of the resonance defined by \eqref{ResSize}. \label{ResonanceComp}}
\end{figure}

Next, we plot the resonance for a slightly longer continuation simulation in
Figure~\ref{ResonanceComp} and compare it for a rational ratio of $\frac{\sigma}{\kappa}$
and an irrational one. From this first investigation it seems that the resonances
for a rational ratio of $\frac{\sigma}{\kappa}$ are evenly spaced for increasing
value of $T$, but for the irrational ratio, the resonances appear to occur more randomly.

When $\frac{\sigma}{\kappa}=m/n$, $\sigma$ and $\kappa$ are both integer multiples of $\kappa/n$.
We would thus expect to see resonances when $T/(\kappa/n)$ is rational, and strong resonances when it is an integer, which has been confirmed experimentally \cite{FleurThesis}. 
However, it appears that not {\em every} integer value of $T/(\kappa/n)$
produces an equally strong resonance.

\part{Backward error analysis} \label{sec:bea}
To better understand how well multisymplectic discretisations preserve 
travelling wave solutions for nonlinear equations, we turn to backward error analysis.
Indeed, similar studies have been performed in the past.
Specifically, implicit midpoint time discretisations of the KdV equation \cite{DFSS} and the
nonlinear {S}chr\"{o}dinger equation \cite{DSS} initialized with an exact travelling wave 
have better properties of error propagation than other (non-conservative) methods.
From another point of view, initializing (\ref{eq:5pt}) with the heteroclinic travelling wave
solution of the modified equation for the scheme 
produces a numerical solution that is accurate to higher order \cite{41MSI}.    
In the following subsection we consider modified equations for the same discretisation (\ref{eq:5pt})
and show that they have Hamiltonian structure with periodic travelling wave solutions that
match the computed travelling wave to higher order.  Then, we show how this idea can be 
extended to other multi-symplectic discretisations and other multi-Hamiltonian PDEs.

\section{Smooth but arbitrary $V^{\prime}$} \label{sec:multistep}

We regard the coincidence between \eqref{DTWE} and a symmetric linear multistep method, 
that holds when $\sigma/\kappa$ is rational, merely as an {\em analogy}. Since backward error analysis explains a good part of 
long-time behaviour of symmetric linear multistep methods (in particular, of their underlying one-step method), 
we apply backward error analysis directly to \eqref{DTWE} for arbitrary $\sigma/\kappa$.
We first write
the second order ODE \eqref{eq:reduced} as
\begin{equation}
\label{f(y)NLW}
\ddot{y}=f(y), \qquad f(y)=\frac{-V'(y)}{c^2-1},
\end{equation}
and the discrete travelling wave equation \eqref{DTWE} as
\begin{equation}
\label{DTWEf(y)}
\frac{\frac{c^{2}}{\kappa^2}\left(y(t+\kappa)-2y(t)+y(t-\kappa)\right)
-\frac{1}{\sigma^2}\left(y(t+\sigma)-2y(t)+y(t-\sigma)\right)}{c^2-1} = f(y(t)).
\end{equation}

Using Taylor series expansions, we get
\[ y(t \pm \kappa) = e^{\pm \kappa D}y(t) \qquad \textup{and} \qquad y(t\pm \sigma) = e^{\pm \sigma D}y(t), \]
where $D$ is the total derivative with respect to $t$.
Substituting these into \eqref{DTWEf(y)} we get
\begin{equation} \label{LinearOp}
L_{\kappa,\sigma}y(t) :=
\frac{\frac{c^{2}}{\kappa^{2}}\left(e^{\kappa D}-2+e^{-\kappa D}\right)-\frac{1}{\sigma^{2}}\left(e^{\sigma D} -2+e^{-\sigma D}\right)}{c^2-1}y(t)=f(y(t))
\end{equation}
were we define the linear operator $L_{\kappa,\sigma}$.
Thus, there exist functions $f_{i}(\varphi,\dot{\varphi})$ with $(i-1)$th
order coefficients of $\sigma$ and $\kappa$, such that
every solution of the modified differential equation
\begin{equation}
\label{SymMDE}
\ddot{y}=f(y)+h^{2}f_{3}(y,\dot{y})+h^{4}f_{5}(y,\dot{y})+\ldots
\end{equation}
satisfies (\ref{DTWE})
up to $O(h^p)$ for all $p$, where $f(y)$ is given in \eqref{f(y)NLW},
and $h=1$ is used as a placeholder for separating terms of different orders.
%
%
To illustrate, we carry out the required calculations up to the 4th order term
by expanding the linear operator $L_{\kappa,\sigma}$ in a Taylor series up to the $4$th power of $h$.
Rearranging \eqref{LinearOp} and applying the derivative $D$ twice gives
\begin{equation}
\label{ExDTWE}
\ddot{y}= D^{2} L_{\kappa,\sigma}^{-1} f(y) = (1+\mu_{2}(\kappa,\sigma)h^{2}D^{2}+\mu_{4}(\kappa,\sigma)h^{4}D^{4}+\dots)f(y).
\end{equation}
where
\begin{equation*}
\mu_{2} = \frac{\sigma^{2}-c^{2}\kappa^{2}}{12(c^{2}-1)},\qquad
\mu_{4} = \frac{3c^{2}\kappa^{4}+3\sigma^{4}+2c^{2}(\kappa^{4}-5\kappa^{2}\sigma^{2}+\sigma^{4})}{720(c^{2}-1)^{2}}.
\end{equation*}
Here the $h^2$ and $h^4$ terms are used as placeholders for ease of calculation and comparison.
We will set $h=1$ at the end of the comparison to the modified differential equation \eqref{SymMDE}.

Now, the derivatives $D^{2}f(y)$ and $D^{4}f(y)$ need to be calculated,
giving expressions involving $f(y)$ and its derivatives and also derivatives of $y$, i.e.
\eqref{SymMDE} for $\ddot{y}$, giving
\begin{align*}
D^{2}f(y) &= f''(y)(\dot{y},\dot{y})+f'(y)\ddot{y} 
=f''(y)(\dot{y},\dot{y})+f'(y)f(y)+ {\cal O}(h^{2}).
\end{align*}
Substituting this back into \eqref{ExDTWE} 
and matching the coefficients of $h^{2}$ with \eqref{SymMDE}, we get
\begin{align}
\label{f3eq}
f_{3}(y,\dot{y}) 
&= \mu_{2}(f''(y)(\dot{y},\dot{y})+f'(y)f(y)).
\end{align}
At order $h^4$ we find $D^{4}f(y)$ and comparing the coefficients of $h^4$ gives the result
\begin{align*}
f_{5}(y,\dot{y})=\mu_{2}(f'(y)f_{3}(y,\dot{y}))&+\mu_{4}(f'(y)^{2}f(y)+3f''(y)f(y)^{2}+5f''(y)(\dot{y},\dot{y})f'(y) \\ &+6f^{(3)}(y)(\dot{y},\dot{y})f(y)+f^{(4)}(y)(\dot{y},\dot{y},\dot{y},\dot{y})).
\end{align*}
Each $f_n$ can be found similarly.

Now, it can be checked through direct calculation that
the modified differential equation \eqref{SymMDE} for the discrete travelling wave equation can be
written up to any power of $h$ as a first-order noncanonical Hamiltonian system
\begin{equation}
\label{eq:modJH}
\dot{z}=\widetilde{J}^{-1}({z})\nabla\widetilde{H}({z}),
\end{equation}
where $z\in\mathbb{R}^{2N}$, $\widetilde J\in\mathbb{R}^{2N\times 2N}$ is a symplectic structure matrix, and
$\widetilde{H}(z)$ is the modified Hamiltonian with expansion
\begin{equation}
\label{eq:modH}
\widetilde{H}(y,\dot{y}) = H(y,\dot{y})+h^{2}H_{3}(y,\dot{y})+h^{4}H_{5}(y,\dot{y})+\ldots.
\end{equation}
For $N=1$, the modified Hamiltonian has
\begin{align*}
H_3 &= \frac{1}{2}\mu_{2}f^2(y)-\mu_{2}f'(y)p^2 \\
H_5 &= \mu_{4}f^2(y)f'(y)+\frac{1}{2}(\mu_{2}^{2}-3\mu_{4})f'(y)^2p^{2}-2\mu_{4}f(y)f''(y)\dot{y}^2-\mu_{4}p^4f^{(3)}(y),
\end{align*}
where $z=(y,p)$, and the modified symplectic structure is
\begin{equation}
\label{eq:modJ}
J^{-1} =
\begin{bmatrix}
0 & \widetilde{K} \\
-\widetilde{K} & 0
\end{bmatrix},
\end{equation}
where $\widetilde{K} =1+h^2K_{3}+h^4K_{5}+\cdots$, and
\[
K_3 = 2 \mu_2 f'(y), \qquad
K_5 = -(\mu_2^2-3\mu_4)f'(y)^2 + 4 \mu_4 f(y)f''(y) - 4 \mu_4 p^2 f'''(y).
\]

Thus the travelling waves of the modified equation can be read off from the level sets of $\widetilde H$.
As this is a smooth function close to $H$, we see that, in the sense of backward error analysis,
all periodic and all heteroclinic travelling waves are preserved by the integrator.

\begin{figure}[h]
\centering
\vspace{-20mm}
\includegraphics[angle=270,width=1\textwidth]{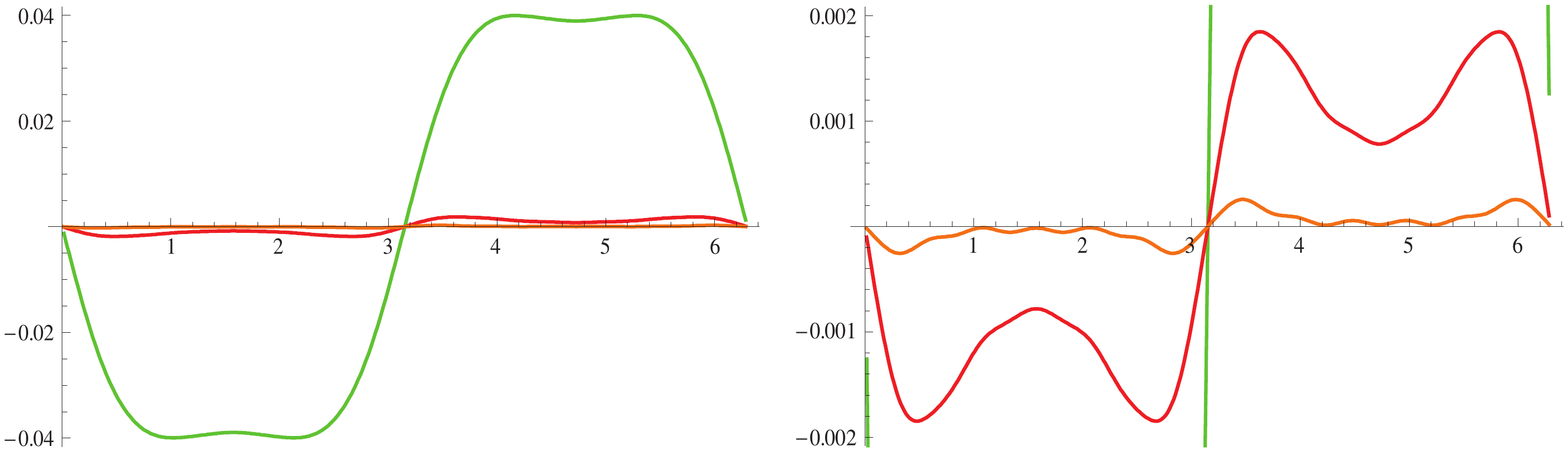}
\vspace{-30mm}
\caption{The errors $y-y_1$ (green), $y-y_2$ (red), and $y-y_3$ (orange) for a 
 non-resonant periodic travelling wave of the form depicted in Figure \ref{Continuation2} 
 (zoomed on right).\label{TruncatedSolutions}}
\end{figure}

We can get some idea of the validity of the method by comparing a computed travelling
wave $y$ with the solutions of the modified differential equation determined by backward
error analysis.
Figure~\ref{TruncatedSolutions} gives these solutions for
$T=2\pi$, $\sigma=\kappa=0.5$, where
\begin{align*}
y_{1} \text{ is the solution of } \ddot{y}&=f(y) \\
y_{2} \text{ is the solution of } \ddot{y}&=f(y)+\mu_{2}f_{3}(y,\dot{y}) \\
y_{3} \text{ is the solution of } \ddot{y}&=f(y)+\mu_{2}f_{3}(y,\dot{y})+\mu_{4}f_{5}(y,\dot{y}),
\end{align*}
which are computed numerically using NDSolve in Mathematica, with $f(y) = \sin(y)/(1-c^2)$. Although the space
and time steps are fairly large, the modified equation gives an excellent approximation to the
actual travelling waves, with more terms improving the accuracy.

\begin{figure}[h]
\centering
\includegraphics[width=0.49\textwidth]{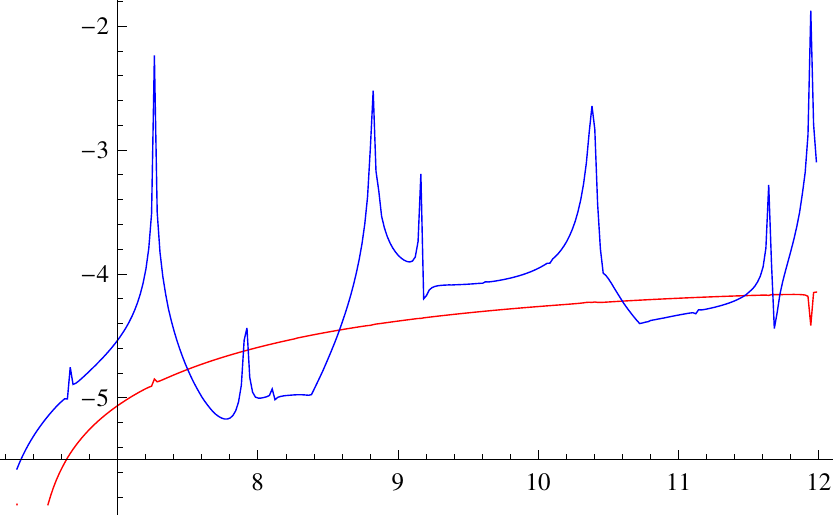}
\includegraphics[width=0.49\textwidth]{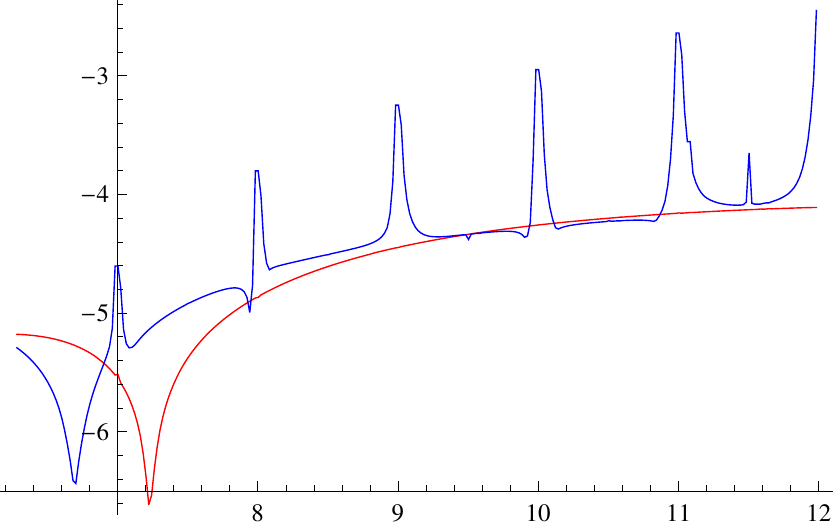}
\caption{Resonance measure $R$ of the numerical solution from Section~\ref{sec:smooth} given by the blue curve.
Resonance measure $R$ of the modified differential equation for the discrete travelling wave equation
given by the red curve. \label{ResonanceAndMDE}}
\end{figure}

We next compare, in Figure \ref{ResonanceAndMDE},
the resonance measure $R$ given in equation \eqref{ResSize}
of the actual travelling waves and those of the modified equation.
We see that the modified equation does not capture the
resonant behaviour. (The same phenomenon is seen for ODEs).
However, the modified equation does approximate the actual travelling waves
very accurately in this sensitive measure.

\begin{figure}[h]
\centering
\includegraphics[width=0.8\textwidth]{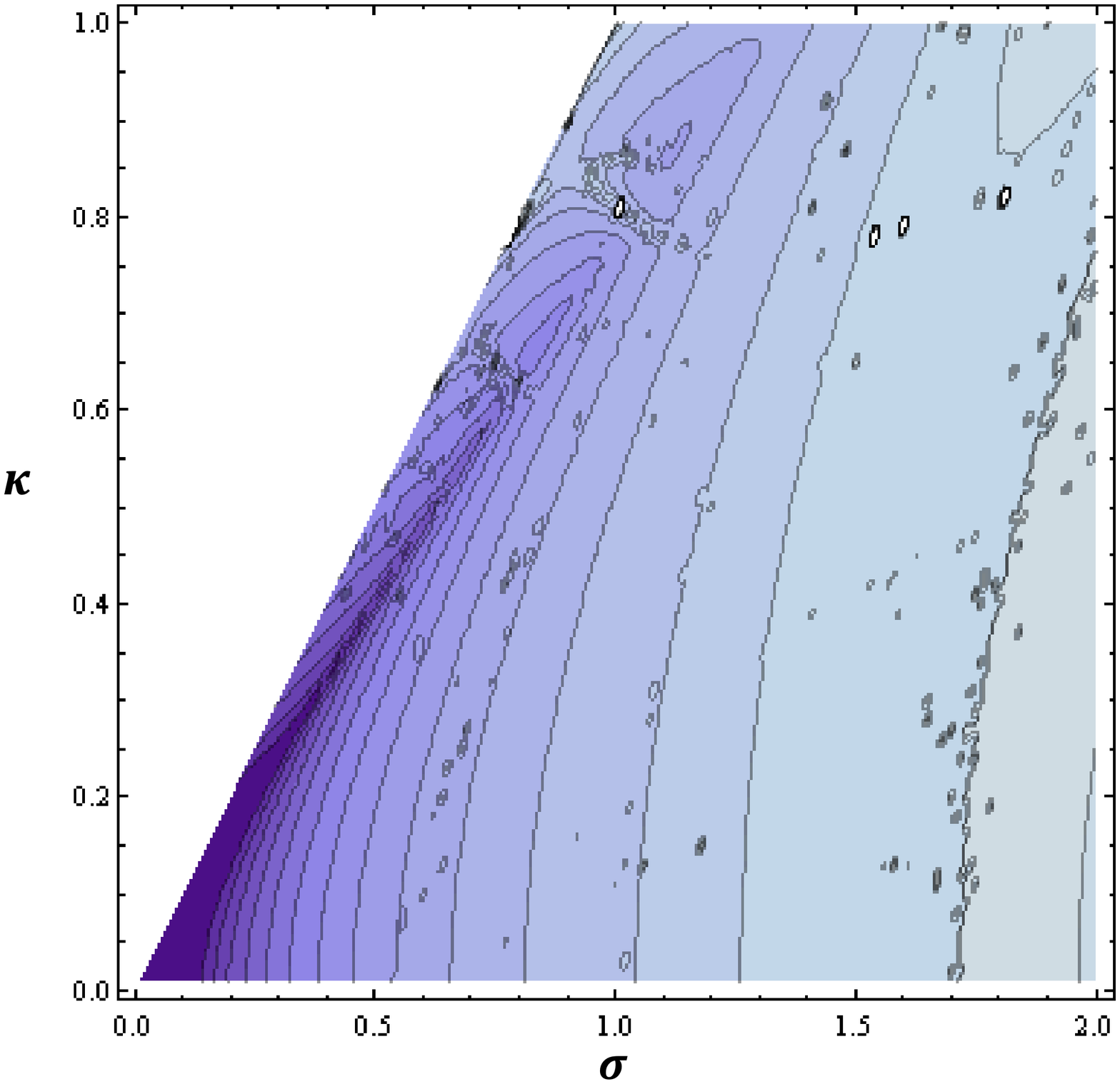}
\caption{Contour plot of $\log \| y - y_2 \|_2$ where $y$ is a discrete travelling wave and
$y_2$ is the corresponding solution of the first truncation of the modified equation. Here $T=2\pi$ and $c=1.3$.\label{ErrorGraphs}}
\end{figure}

Finally, for fixed values of $T$ and $c$, we compare the actual travelling waves and those
predicted from the modified equation for all $\sigma$ and $\kappa$. Figure \ref{ErrorGraphs}
shows a contour plot of $\log \| y - y_2 \|_2$ where $y$ is a discrete travelling wave and
$y_2$ is the corresponding solution of the first truncation of the modified equation.
Overall, the error behaviour is that predicted by the Taylor series (the smooth part of the figure), and
the resonances show up as small anomalies. Only two such anomalies (near $(\sigma,\kappa)=(0.8,0.6)$ and $(1.0,0.8)$) are at all significant.

\section{Generalizations to other methods and PDEs} \label{sec:conclude}

The same method of backward error analysis can be used to approximate the travelling wave solutions
of any multisymplectic integrator for any multi-Hamiltonian problem.
We illustrate this using the Preissman box scheme, a multisymplectic integrator that can
(unlike the leapfrog method considered earlier) be applied to any multi-Hamiltonian PDE.
This gives
\begin{align}
\label{MHPDEDis}
(K D_{t}M_{x}+L D_{x}M_{t})z_{i}^{n}&=\nabla S(M_{t}M_{x}z_{i}^{n}),
\end{align}
where $D_t z_i^n = (z_i^{n+1}-z_i^n)/(\Delta t)$, $D_x z_i^n = (z_{i+1}^n - z_i^n)/(\Delta x)$, $M_x z_i^n = (z_i^n + z_{i+1}^n)/2$, and $M_t = (z_i^n + z_i^{n+1})/2$.

Substituting travelling wave coordinates  we get the difference operators
\begin{align*}
D_{\xi}^{\kappa}\varphi(\xi)=-\frac{c(\varphi(\xi)-\varphi(\xi-\kappa))}{\kappa}, \qquad D_{\xi}^{\sigma}\varphi(\xi)=\frac{\varphi(\xi+\sigma)-\varphi(\xi)}{\sigma}, \\
M_{\xi}^{\kappa}\varphi(\xi)=\frac{(\varphi(\xi)+\varphi(\xi-\kappa))}{2}, \qquad M_{\xi}^{\sigma}\varphi(\xi)=\frac{\varphi(\xi)+\varphi(\xi+\sigma)}{\sigma},
\end{align*}
so that upon substituting travelling wave coordinates in the discretisation \eqref{MHPDEDis} we get the system of difference equations,
\begin{equation}
\label{eq:3}
(K D_{\xi}^{\kappa}M_{\xi}^{\sigma}+L D_{\xi}^{\sigma}M_{\xi}^{\kappa})\varphi(\xi)=\nabla S(M_{\xi}^{\kappa}M_{\xi}^{\sigma}\varphi(\xi)).
\end{equation}
If the linear operator on the left is  nonsingular, we can move it to the right hand side and apply the $\xi$ derivative $D$ to both sides to get
\begin{equation}
\label{eq:gen}
\dot\varphi(\xi) = D((K D_{\xi}^{\kappa}M_{\xi}^{\sigma}+L D_{\xi}^{\sigma}M_{\xi}^{\kappa})^{-1}\nabla S(M_{\xi}^{\kappa}M_{\xi}^{\sigma}\varphi(\xi)).
\end{equation}
As before, the right hand side can be expanded in a Taylor series in $h$
(or $\sigma$, $\kappa$) which allows the modified equation to be determined term-by-term.
We expect that it will be conjugate to symplectic and thus the preservation of travelling waves,
in the sense of backward error analysis, can be determined from the Hamiltonian of the continuous travelling wave equation.

This occurs, for example, in the nonlinear Schr\"odinger equation
\begin{equation*}
i\psi_{t}+\psi_{xx}+2\left|\psi\right|^{2}\psi=0, \qquad \psi\in\mathbb{C}.
\end{equation*}
Setting $\psi = p + i q$, $\psi_x = v + i w$,
the multi-Hamiltonian form is (\ref{eq:PDE}) with
\begin{equation*}
z=
\begin{bmatrix}
p \\
q \\
v \\
w
\end{bmatrix}, \qquad
K =
\begin{bmatrix}
0 & -1 & 0 & 0 \\
1 & 0 & 0 & 0 \\
0 & 0 & 0 & 0 \\
0 & 0 & 0 & 0
\end{bmatrix}, \qquad
L =
\begin{bmatrix}
0 & 0 & 1 & 0 \\
0 & 0 & 0 & 1 \\
-1 & 0 & 0 & 0 \\
0 & -1 & 0 & 0
\end{bmatrix}
\end{equation*}
and $S=-\frac{1}{2}(p^2+q^2)-\frac{1}{2}(v^2+w^2)$.
Equation \eqref{eq:gen} in this case becomes
\begin{equation*}
\dot{\varphi}(\xi)=D
\begin{bmatrix}
0 & 0 & -\frac{1}{D_{\xi}^{\sigma}M_{\xi}^{\kappa}} & 0 \\
 & 0 & 0 & -\frac{1}{D_{\xi}^{\sigma}M_{\xi}^{\kappa}} \\
\frac{1}{D_{\xi}^{\sigma}M_{\xi}^{\kappa}} & 0 & 0 &  -\frac{D_{\xi}^{\kappa}M_{\xi}^{\sigma}}{(D_{\xi}^{\kappa}M_{\xi}^{\sigma})^{2}} \\
0 & \frac{1}{D_{\xi}^{\sigma}M_{\xi}^{\kappa}} & \frac{D_{\xi}^{\kappa}M_{\xi}^{\sigma}}{(D_{\xi}^{\kappa}M_{\xi}^{\sigma})^{2}} & 0
\end{bmatrix}\nabla S(M_{\xi}^{\kappa}M_{\xi}^{\sigma}\varphi(\xi)).
\end{equation*}

The modified equation of the discrete travelling waves is a 4-dimensional noncanonical Hamiltonian system.
Both the continuous and the modified systems are rotationally invariant.
Thus, travelling waves are generically preserved in this example.
For 4-dimensional systems without such a symmetry, more complicated situations can arise.
For example, there can be situations in which the continuous travelling wave equation is integrable while the modified
equation for discrete travelling waves is not integrable. In such cases not all travelling waves would be preserved.
However, even in these cases the backward error analysis gives considerable insight into the preservation of dynamics.

The method can also be applied if the linear operator in \eqref{eq:3} is singular. We illustrate this
for the nonlinear wave equation (\ref{eq:PDE}) in the formulation
\begin{equation*}
z =
\begin{bmatrix}
u \\
v \\
w
\end{bmatrix}, \qquad
K  =
\begin{bmatrix}
0 & 1 & 0 \\
-1 & 0 & 0 \\
0 & 0 & 0
\end{bmatrix}, \qquad
L  =
\begin{bmatrix}
0 & 0 & -1 \\
0 &  0 & 0 \\
1 & 0 & 0
\end{bmatrix},
\end{equation*}
and $S(z) = -V(u) + \frac{1}{2}(w^{2} - v^{2})$.

The multi-Hamiltonian discretisation \eqref{MHPDEDis} for the leapfrog discretisation of the nonlinear wave equation may be written
\begin{equation}
\label{MHPDELeapDis}
(K D_{t}+L D_{x})z_{i}^{n}=\nabla S(z_{i}^{n})
\end{equation}
where
\begin{equation}
\label{1stOrderLF}
D_{t}z_{i}^{n}=\frac{z_{i}^{n+\frac{1}{2}}-z_{i}^{n-\frac{1}{2}}}{\Delta t}, \qquad D_{x}z_{i}^{n}=\frac{z_{i+\frac{1}{2}}^{n}-z_{i-\frac{1}{2}}^{n}}{\Delta x}.
\end{equation}
In matrix form \eqref{MHPDELeapDis} becomes
\begin{equation*}
\begin{bmatrix}
0 & D_{t} & -D_{x} \\
-D_{t} & 0 & 0 \\
D_{x} & 0 & 0
\end{bmatrix}z_{i}^{n}=\nabla S(z_{i}^{n}).
\end{equation*}
Letting
\begin{equation}
\label{GenMHTWC}
D_{\xi}^{\gamma}\varphi(\xi)=\frac{\varphi(\xi+\frac{1}{2}\gamma)-\varphi(\xi-\frac{1}{2}\gamma)}{\gamma},
\end{equation}
the discrete travelling wave equation is
\begin{equation*}
(L D_{\xi}^{\sigma}-cK D_{\xi}^{\kappa})\varphi(\xi)=\nabla S(\varphi(\xi)),
\end{equation*}
and in matrix form,
\begin{equation}
\label{DisMHPDE}
\begin{bmatrix}
0 & -cD_{\xi}^{\kappa} & -D_{\xi}^{\sigma} \\
cD_{\xi}^{\kappa} & 0 & 0 \\
D_{\xi}^{\sigma} & 0 & 0
\end{bmatrix}\varphi(\xi)=\nabla S(\varphi(\xi)),
\end{equation}

The operator on the left hand side is singular. We eliminate $\varphi_3$ using the last
equation, move the remaining linear operator to the right hand side, and apply $D$ to both sides to get
\begin{equation*}
\dot{\varphi}=D
\begin{bmatrix}
0 & (cD_{\xi}^{\kappa})^{-1} \\
-cD_{\xi}^{\kappa}(c^2(D_{\xi}^{\kappa})^2-(D_{\xi}^{\sigma})^2)^{-1} & 0
\end{bmatrix}
\begin{bmatrix}
-V'(\varphi_{1}) \\
-\varphi_{2}
\end{bmatrix}.
\end{equation*}
The  right hand side can be expanded in a Taylor series to produce 
the modified equation, which is used to study the phase portraits of the modified Hamiltonian.

It can be seen that the ideas and methods developed here  apply
to higher-dimensional wave equations and to the preservation of travelling waves of general multi-Hamiltonian systems by very large classes of multisymplectic integrators.

\section{Discussion}
The standard application of backward error analysis for symplectic integrators
reduces an $N$-dimensional symplectic map to an $N$-dimensional Hamiltonian flow.
This is essentially a reduction by one dimension. The present application of
backward error analysis, in conjunction with the travelling wave reduction,
has reduced an infinite-dimensional map to a two-dimensional
Hamiltonian ODE, a far more spectacular reduction. This is even more striking in view
of the fact that there is no accepted backward error analysis for partial difference equations in general
(see \cite{63BEA}).

Thus, we see that backward error analysis is exceptionally powerful for this problem,
reducing the intractable functional equation \eqref{DTWE} to a simple
planar Hamiltonian ODE. Furthermore, by comparing the phase portraits of the original
and modified systems we see that any orbits that are structurally stable,
that is any orbits that are also present up to small perturbations in nearby Hamiltonian systems,
are preserved. This includes periodic and heteroclinic travelling waves for the sine-Gordon equation,
 although this result is not specific to the sine-Gordon equation or its
integrability and holds for  generic potentials. This preservation holds in the sense of
backward error analysis, that is, up to any power of the time and space step sizes.

This situation in the ODE case, analyzing the preservation of phase portraits by
symplectic integrators, can be studied  using the language of topological equivalence.
It is shown in \cite{MPQ} that for generic smooth Hamiltonians $H$, the level sets of $H$ and 
of the modified Hamiltonian $\widetilde H$ of a symplectic integrator are topologically equivalent.
It is in this sense that, in the two-dimensional case, symplectic integrators preserve all orbits.
Thus, we can expect, for example, that if a travelling wave corresponds
to a topologically unstable orbit of the travelling wave equation,
then it will not be preserved by typical multisymplectic integrators.
This, however, is an exceptional case, and even in higher dimensions we
can expect that typical travelling waves, whose existence is due to preserved
features of the travelling wave equation like dimension and linear symmetries,
will be preserved by multisymplectic integrators.

Although the numerical examples given above suggest that the modified vector field
is a very good approximation to the travelling waves, which in direct numerical
calculations do appear to exist, of course this is no substitute for a rigorous proof
of their existence. Instead, we regard the backward error analysis as suggesting the
mechanism by which travelling waves can be preserved and as a guide to the development
of good numerical methods.

\vspace{4ex}
\noindent
{\Large \bf Appendix: Steady state solutions}\\

Parts I and II of this paper have exclusively considered multisymplectic methods. In this appendix
we consider non-symplectic spatial discretizations. These necessarily lead to non-multisymplectic full discretizations.
The steady state solutions of both types of methods can be understood and described in great detail, which
allows a comparison of the two types. 

First, consider the steady state solutions of the 5-point discretization (\ref{eq:5pt}).
These correspond to travelling waves of zero speed. 
They obey the equation
\begin{equation}
\label{eq:5ptss}
-\frac{1}{(\Delta x)^2}\left(u_{i+1}^{n}-2u_{i}^{n}+u_{i-1}^{n}\right) = - V'(u_{i}^{n}),
\end{equation}
which is just the leapfrog method applied to $-u_{xx} = -V'(u)$
whose behaviour is well known. We summarize it here from two points of view.

First, for small $\Delta x$, backward error analysis shows how solutions to (\ref{eq:5ptss})
satisfy extremely closely (up to terms exponentially small in $\Delta x$) a modified equation
which is also a planar Hamiltonian system. Thus, in this sense of backward error analysis,
comparing the steady-state solutions of the 5-point method with those of the 
PDE $u_{tt} = u_{xx} - V'(u)$ amounts to comparing the phase portraits of two planar
Hamiltonian systems with nearby Hamiltonians. Generic periodic orbits and hetero/homoclinic orbits persist.

Second, for finite $\Delta x$, the phase portrait of the leapfrog method for generic potentials is,
to some extent, well understood. Where the continuous problem has a hetero/homoclinic orbit, the discrete
problem will have two orbits. Likewise, periodic orbits of the discrete problem come in pairs that
approximate continuous periodic orbits of the continuous problem. These provide periodic solutions of
the discrete problem on infinite and on periodic domains. According to KAM theory, strongly nonresonant
periodic solutions of the continuous problem persist and thus provide periodic discrete solutions on an
infinite domain.

In contrast, non-symplectic spatial discretizations, such as upwind differences,
show a completely different behaviour. Their steady state equation corresponds to
the application of a non-symplectic integrator, such as Euler's method or a Backwards
Differentiation Formula (BDF), to $-u_{xx}=-V'(u)$.
Typically, all hetero/homoclinic orbits, and all but a finite number
of periodic orbits (typically one), are destroyed by such a discretization.

We can also use the steady-state case to compare multisymplectic schemes with {\em symmetric} ones.
The relationship between symplectic and symmetric methods is now well understood in numerical ODEs,
but is not always appreciated in numerical PDEs. For example, Cohen \cite[p. 34]{cohen2002} writes
``All the schemes described below will be centered, since uncentered schemes generate numerical
dissipation for wave equations which satisfy a principle of energy conservation.''
In fact, symmetric (i.e. ``centered'') methods only preserve symmetric waves and do
not typically preserve nonsymmetric waves.



\begin{figure}
\centering
\includegraphics[scale = 0.23]{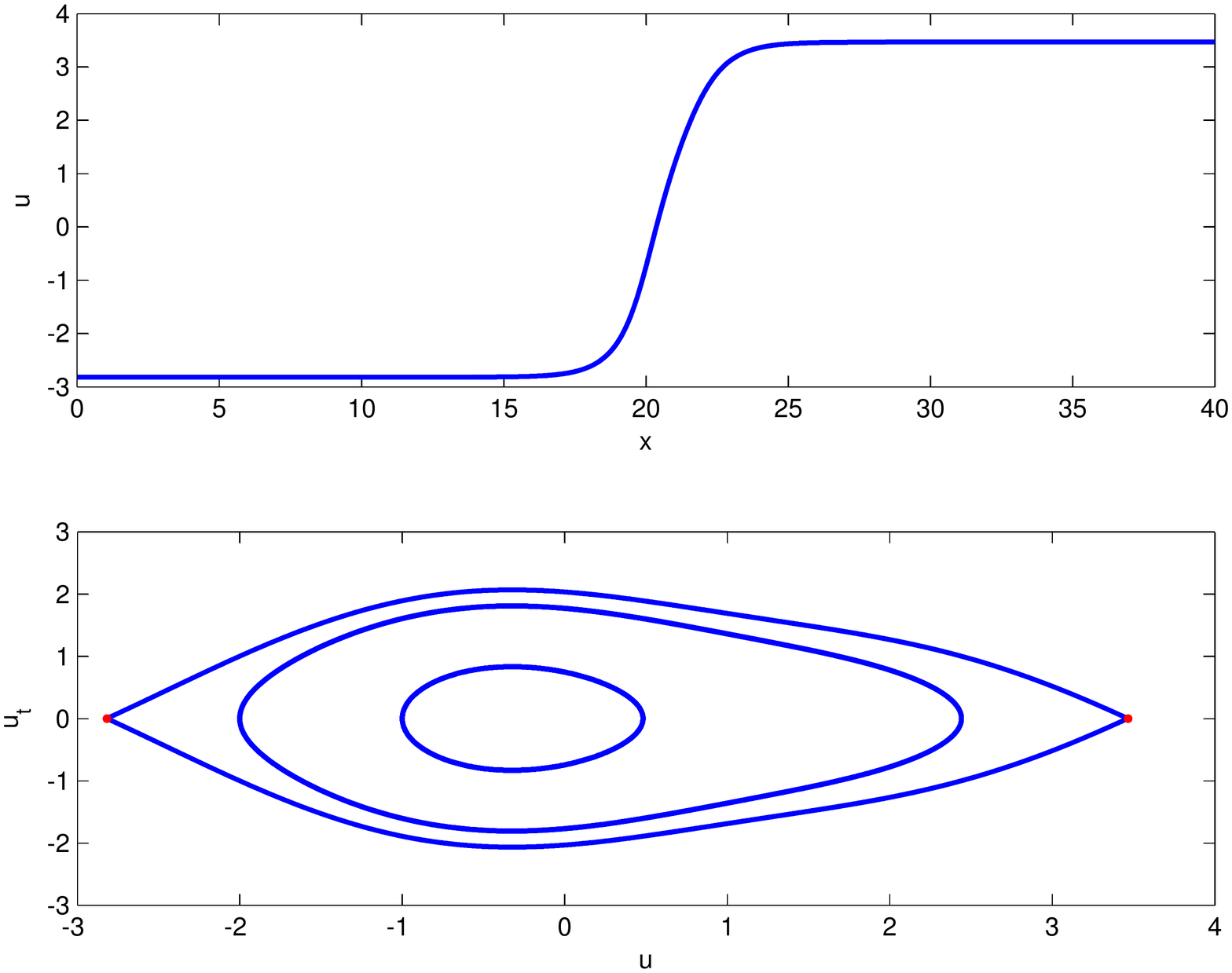}
\includegraphics[scale = 0.23]{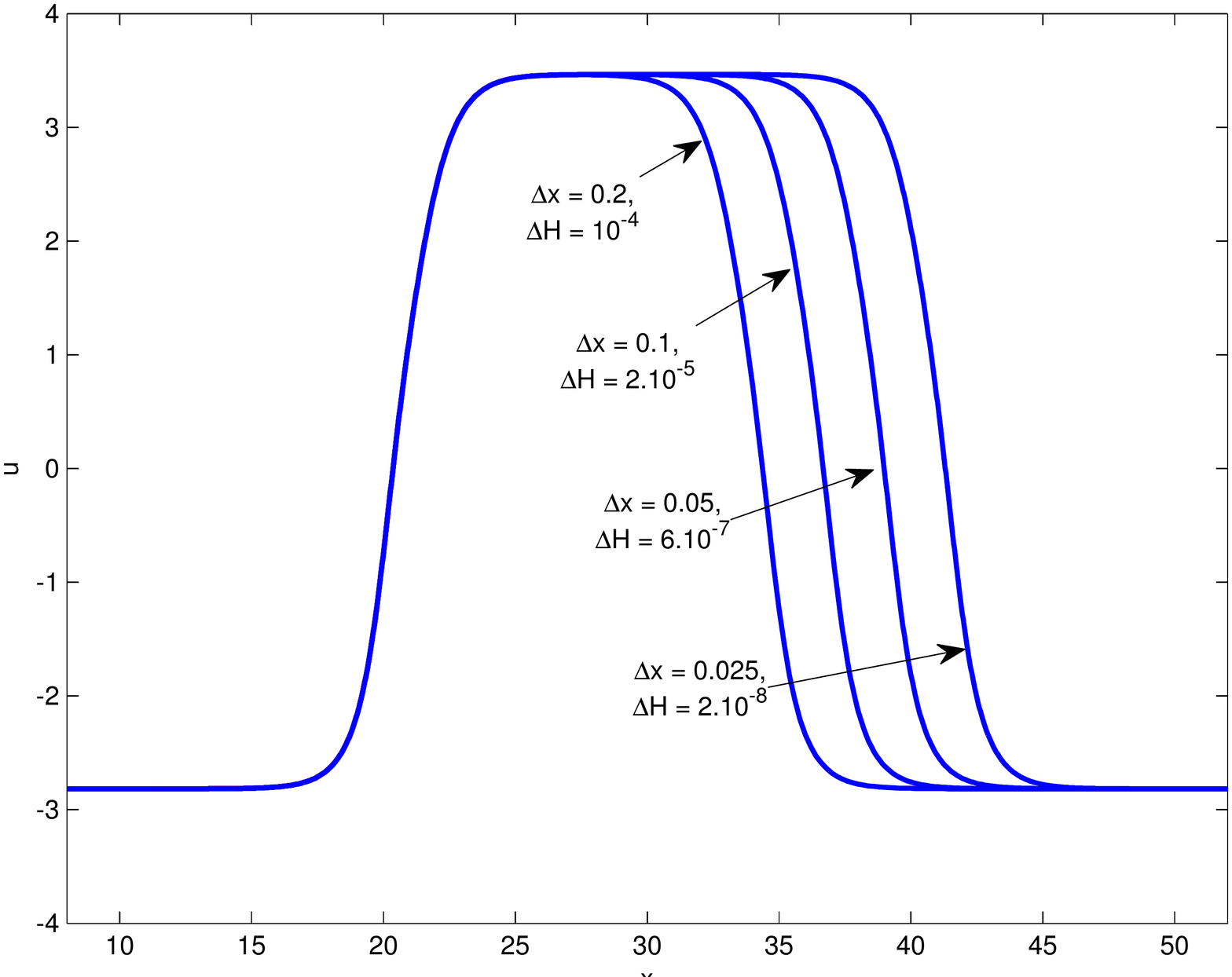}
\caption{Left: Exact solutions of steady-state solutions of \eqref{eq:NLW} with $-V'(u)=\sin u + \frac{2}{5}\cos 2u$ 
(top: one of two heteroclinic solutions; bottom: phase portrait showing two heteroclinic and two periodic orbits.)
Right: Steady state solutions of a symmetric spatial discretisation (Lobatto IIIA); there are homoclinic, but no heteroclinic orbits for any $\Delta x$.\label{fig:SymSol}}
\end{figure}

Consider the force $-V'(u) =\sin u + \frac{2}{5}\cos2u$ with $u\in\mathbb{R}$, 
used in
 \cite{8GNI} to study non-conservation of energy by symmetric integrators. 
The phase portrait of the steady-state equation $-u_{xx} = -V'(u)$  is shown in Figure~\ref{fig:SymSol}.
From this we see that the PDE has both periodic  and heteroclinic steady states, but no homoclinic steady states.
One of the heteroclinic steady states is shown on the left-hand side of Figure~\ref{fig:SymSol}.
The ODE $-u_{xx} = -V'(u)$ does have a reversing symmetry, namely $(u,x,t)\mapsto (u,-x,-t)$,
but its heteroclinic orbits are not invariant under this reversing symmetry. That is, they
are nonsymmetric orbits.

As an example of a spatial discretization which is symmetric but not symplectic, we will use 
the 3-stage Lobatto IIIA method (a  Runge--Kutta method).
We calculate the solutions of this method that are asymptotic to the fixed point near $u=-\pi$
and find that the method has no heteroclinic steady states for any of the values of $\Delta x$
considered, and that it does have homoclinic steady states. These are shown on the right-hand side of Figure~\ref{fig:SymSol}. That is, this spatial discretization has a
qualitatively wrong steady state structure for all $\Delta x$.

These considerations demonstrate the the superiority of multisymplectic over non-multisymplectic methods in the preservation of qualitative solution structure.


\end{document}